\documentclass[12pt]{article} 
\usepackage{amsfonts}           
\usepackage{myart98}
\title{Multiparametric Dissipative Linear Stationary Dynamical Scattering 
Systems: Discrete Case}
\author{Dmitriy S. Kalyuzhniy}
\date{}
\newcommand{\nspace}[2]{\ensuremath{{\mathbb{#1}}^{#2}}}
\newcommand{\Hspace}[1]{\ensuremath{\mathcal{#1}}}
\def\ifundefined#1{\expandafter\ifx\csname#1\endcsname\relax}
\newcommand{\comment}[1]{}

\usepackage{amsfonts}

\newcommand{\Cliff}[2][\comment]{\ensuremath{{\cal C}\kern-
0.18em\ell(#1,#2)}}

\ifundefined{qed}
    \DeclareMathSymbol{\qed}{0}{AMSa}{"03}
\fi
\newcommand{\scalar}[2]{\langle #1,#2\rangle}
\newcommand{\bigscalar}[2]{\left\langle #1,#2\right\rangle}

\providecommand{\eqref}[1]{\textup{(\ref{#1})}}

\begin{document}
\maketitle

\begin{abstract}
We propose the new generalization of linear stationary dynamical systems
with discrete time $t\in\mathbb{Z}$ to the case $t\in\nspace{Z}{N}$. The
dynamics of such a system can be reproduced by means of its associated
multiparametric Lax-Phillips semigroup. We define multiparametric
dissipative, and conservative scattering systems and interpret them in terms
 of
operator
colligations, of the associated semigroup and of ``energy'' relations
for
system data. We prove the Agler's type theorem describing the class of
holomorphic operator-valued functions on the polydisc $\nspace{D}{N}$
that
are the transfer functions of multiparametric conservative scattering
systems.
\begin{description}
 \item[Keywords:] \emph{Dissipative systems, multiparametric 
Lax-Phillips
semigroup, generalized Schur class, conservative realizations}
 \item[AMS Subject Classification:] 93C35, 93B15, 32A10, 47A48.
\end{description}
\end{abstract}

\setcounter{section}{-1}
\section{Introduction and preliminaries}
In the present paper we introduce the concept of multiparametric linear
stationary dynamical system (LSDS) and, in particular, the concept of
dissipative (conservative) scattering LSDS, establish some properties of
dissipative
systems, give the description of the class of contractive
operator-valued
functions holomorphic on the open unit polydisc $\nspace{D}{N}$ which
are the
transfer functions of multiparametric conservative scattering LSDSs. We
are
based on the results of J.~Agler on the characterization of the
generalized
Schur
class (see \cite{Ag}).

In our considerations the multidimensional  parameter
$t\in\nspace{Z}{N}$
plays a role of
``multidimensional time'' since the introduced concepts
generalize the concept of LSDS with (one-dimensional) discrete time and,
in
particular, the concept of dissipative (conservative) scattering LSDS. In this
introductory section we shall recall the basic
definitions and the most important facts of the theory of one-parametric
dissipative scattering systems (see \cite{A2} or \cite{BC}) and give the
motivation of such generalization.

Let us start with the standard definition of  \emph{LSDS} for the case
$t\in\mathbb{Z}$ as the following system of equations:
\begin{equation}\label{eq:sys}
\alpha :\left\{ \begin{array}{lll}
                 x(t+1)&=&Ax(t)+Bu(t),\\
                 y(t)&=&Cx(t)+Du(t),
                \end{array}
        \right.\quad
(t=0,1,2,\ldots)
\end{equation}
where $A:\Hspace{X}\to\Hspace{X}$, $B:\Hspace{N^-}\to\Hspace{X}$,
$C:\Hspace{X}\to\Hspace{N^+}$, $D:\Hspace{N^-}\to\Hspace{N^+}$ are
bounded linear operators in separable Hilbert spaces, and the initial
condition
\begin{equation}\label{eq:init-sys}
x(0)=x_0,
\end{equation}
with prescribed $x_0\in\Hspace{X}$. The spaces $\Hspace{X}$,
$\Hspace{N^-}$,
$\Hspace{N^+}$ are called the \emph{state space}, the \emph{input
space},
and the \emph{output space} respectively for the LSDS $\alpha$. The
system
$\alpha$ is called a \emph{dissipative scattering} LSDS if the system
matrix\begin{equation}\label{eq:matr}
G=\left( \begin{array}{ll}
          A & B \\
          C & D
         \end{array}
  \right)
\in [\Hspace{X}\oplus\Hspace{N^-},\Hspace{X}\oplus\Hspace{N^+}]
\end{equation}
(here $[\Hspace{H}_1,\Hspace{H}_2]$ denotes the space of all bounded
linear
operators acting from a separable Hilbert space $\Hspace{H}_1$ into a
separable Hilbert space $\Hspace{H}_2$) defines a contractive operator,
i.e., $G^*G\leq I_{\Hspace{X}\oplus\Hspace{N^-}}$ in sense of Hermitian
operators. The system $\alpha$ is called a \emph{conservative
scattering}
LSDS if $G$ is unitary,
i.e.,\begin{equation}\label{eq:unitary}
G^*G=I_{\Hspace{X}\oplus\Hspace{N^-}},\quad
GG^*=I_{\Hspace{X}\oplus\Hspace{N^+}}. \end{equation} Note that
dissipativity
(resp. conservativity) conditions have physical sense of the
dissipation
(resp. of the full balance) of energy in a system. We hold the
terminology
of \cite{HM}, \cite{W},
 however some authors (see e.g. \cite{A2}) use the term
``passive system'' instead of the term
``dissipative system''. One may gather the data of a system   using
the
notation $\alpha =(A, B, C, D; \Hspace{X}, \Hspace{N^-}, \Hspace{N^+}).$
In
the case of  dissipative (resp. conservative) scattering LSDS the
aggregate
$\alpha =(A, B, C, D; \Hspace{X}, \Hspace{N^-}, \Hspace{N^+})$ is
called a \emph{contractive} (resp. \emph{unitary}) \emph{operator
colligation}.

The operator-valued function \begin{equation}\label{eq:tf-sys}
\theta_\alpha (z)=D+zC{(I_{\Hspace{X}}-zA)}^{-1}B
\end{equation} which has to be considered on some neighbourhood of $z=0$
in $\mathbb{C}$ is called the \emph{transfer function} of a system
$\alpha $. It is known (see \cite{A2} or \cite{BC}) that the transfer
function $\theta_\alpha$ of an arbitrary dissipative (in particular,
conservative) scattering LSDS $\alpha $ belongs to the \emph{Schur
class} $S(\Hspace{N^-}, \Hspace{N^+})$ consisting of all functions holomorphic
on
the open unit disc $\mathbb{D}$
with values equal to contractive operators from $[\Hspace{N^-}, \Hspace{N^+}]$.
The converse
statement is also true (see \cite{Br}, \cite{Sz.-NF}): for an arbitrary
$\theta\in S(\Hspace{N^-}, \Hspace{N^+})$ there exists its conservative
realization, that is a conservative scattering LSDS (a unitary
colligation)
$\alpha $ of which $\theta$ is the transfer function (the characteristic
function of the corresponding unitary colligation), i.e., $\theta
=\theta_\alpha $. And what is more, the following result, fundamental for the
theory
of characteristic operator-valued functions
 (see \cite{Br}), is valid
(we shall recall it here in the systems theory language by \cite{BC}): for an
arbitrary $\theta\in S(\Hspace{N^-}, \Hspace{N^+})$ there exists a \emph{closely
connected} conservative realization
$\alpha =(A, B, C, D; \Hspace{X}, \Hspace{N^-}, \Hspace{N^+})$, i.e.,
 $\alpha $ is a conservative system and\begin{equation}\label{eq:cc-sys}
\Hspace{X}=\bigvee_{p} p(A, A^*)(B\Hspace{N^-}+C^*\Hspace{N^+})
\end{equation}
(here the symbol ``$\bigvee$'' denotes the closure of the linear span of
some subspaces, $p$ runs the set of all monomials in two noncommuting
variables); such realization is uniquely determined by $\theta$ up to
unitary similarity. The close connectedness of a conservative
system is equivalent to the complete non-unitarity of its (contractive)
\emph{main operator} $A$ (see \cite{BC}). Thus there are close relations
between the theory of dissipative (in particular, conservative) scattering
systems, the theory of contractions by Sz.-Nagy -Foias \cite{Sz.-NF},
the theory of unitary colligations \cite{Br} and the theory of
holomorphic functions of a complex variable, namely of Schur class
functions.

The dynamics of LSDS \eqref{eq:sys}--\eqref{eq:init-sys} can be
described
by means of its
associated Lax-Phillips semigroup in the corresponding scattering scheme
(see \cite{LP}, and also \cite{AA}, \cite{A1}, \cite{N}). The
\emph{abstract Lax-Phillips scattering scheme} for the discrete case
includes some separable Hilbert space $\Hspace{H}$, its subspaces
$\Hspace{D_+}$ (the \emph{outgoing space}) and $\Hspace{D_-}$ (the
\emph{ingoing space}), and some bounded linear operator $W$ in
$\Hspace{H}$ (the \emph{generator} of the \emph{Lax-Phillips semigroup}
$\mathfrak{W}=\{ W^t\mid t=0,1,2,\ldots\} $) such that: (i)
$\Hspace{D_+}\perp\Hspace{D_-}$; (ii)
$W\Hspace{D_+}\subset\Hspace{D_+}$, $W^*\Hspace{D_-}\subset\Hspace{D_-
 }$; (iii) the operators $W|\Hspace{D_+}$ and $W^*|\Hspace{D_-}$
are isometric; (iv) $\bigcap_{n=0}^\infty W^n\Hspace{D_+}=\{ 0\}
=\bigcap_{n=0}^\infty W^{*n}\Hspace{D_-}$. We postpone the definition
of the associated Lax-Phillips semigroup $\mathfrak{W}_\alpha$ for LSDS
$\alpha$ till Section~\ref{sec:rem-sys}, we only notice here that the
semigroup
$\mathfrak{W}_\alpha$ is contractive (resp. unitary) if and only if the
system $\alpha $ is dissipative (resp. conservative).

The powerful instrument in the investigation of properties of systems
and their main operators is the notion of dilation (see e.g. \cite{A2},
\cite{BC}). The LSDS $\widetilde{\alpha}=(\widetilde{A}, \widetilde{B},
\widetilde{C}, D; \widetilde{\Hspace{X}}, \Hspace{N^-}, \Hspace{N^+})$
is called a \emph{dilation of the LSDS} $\alpha =(A, B, C, D;
\Hspace{X},\Hspace{N^-}, \Hspace{N^+})$ if the following conditions are
fulfilled:\begin{eqnarray*}
\widetilde{\Hspace{X}}=\Hspace{D}\oplus\Hspace{X}\oplus\Hspace{D_*},
\quad A=P_\Hspace{X}\widetilde{A}|\Hspace{X},\quad
B=P_\Hspace{X}\widetilde{B},\quad
C=\widetilde{C}|\Hspace{X},\\
\widetilde{A}\Hspace{D}\subset\Hspace{D},\quad
\widetilde{C}\Hspace{D}=\{ 0\} ,\quad
\widetilde{A}^*\Hspace{D_*}\subset\Hspace{D_*},\quad
\widetilde{B}^*\Hspace{D_*}=\{ 0\}
\end{eqnarray*}
(here $P_\Hspace{X}$ denotes the orthogonal projector onto
$\Hspace{X}$). In particular, it means that
$A^t=P_\Hspace{X}\widetilde{A}^t|\Hspace{X}\ (t=1,2,\ldots )$, i.e., the
operator $\widetilde{A}$ is a \emph{dilation of the operator} $A$. The
LSDS $\alpha $ is called \emph{minimal} if it is not a dilation of any
other, different from it, system. The notion of minimal system plays an
important role in control theory; in the case of finite-dimensional
spaces $\Hspace{N^-}$ and $\Hspace{N^+}$ and a rational matrix-valued
function $\theta (\cdot )$ there exists a minimal system of all systems
of which the transfer function is $\theta$, and its state space has the
definite finite dimension which is minimal of all state space dimensions
of such systems (see \cite{AFK}). An arbitrary LSDS $\alpha $ is a
dilation of some minimal LSDS $\alpha _{min}$ with the same transfer
function, moreover if $\alpha $ is a dissipative scattering LSDS then
$\alpha _{min}$ is also dissipative (see \cite{A2}). Together with the
theorem
on the conservative realization of a Schur class function this implies
that any Schur class function has a minimal dissipative realization. On the
other hand (see \cite{A2}), each dissipative scattering LSDS $\alpha $
allows a conservative dilation $\widetilde{\alpha }$, that is the system
analogue of the classical B. Sz.-Nagy theorem (see \cite{Sz.-NF}) on the
existence of a unitary dilation for an arbitrary contraction in a
Hilbert space.

Summarizing all the foregoing let us extract the following aspects of
the theory of dissipative scattering LSDSs with (one-dimensional) discrete
time: (a)~the connection with the theory of holomorphic functions of one
complex variable (namely, of Schur class functions on $\nspace{D}{N}$);
(b)~the connection with operator theory (namely, with the theory of
contractive and unitary operators and operator colligations in Hilbert
spaces); (c)~the connection with the Lax-Phillips scattering theory;
(d)~the theory of dilations of systems as an apparatus in control theory
and in the theory of operators (colligations) in Hilbert spaces.

When constructing the theory of multiparametric dissipative scattering LSDSs
it seems important for the above-mentioned aspects of systems theory to
have meaningful generalizations. In our approach this is realized in the
following way. (a$'$)~By means of the notion of the transfer function of
a system the connection with the theory  of holomorphic functions of
several complex variables (namely, of contractive operator-valued
holomorphic functions on $\nspace{D}{N}$) is established. (b$'$)~If
$\zeta
=(\zeta _1,\ldots ,\zeta _N)$ is a varying point on the unit torus
$\nspace{T}{N}$, $\mathbf{A}=(A_1,\ldots ,A_N)$ and
$\mathbf{G}=(G_1,\ldots ,G_N)$ are $N$-tuples of operators from
$[\Hspace{X}, \Hspace{X}]$ and $[\Hspace{X}\oplus\Hspace{N^-},
\Hspace{X}\oplus\Hspace{N^+}]$ respectively, then the pencil of
operators $\zeta\mathbf{A}:=\sum _{k=1}^N\zeta _kA_k$ is an analogue
of the main operator $A$ of a one-parametric LSDS, and the pencil of
block matrices $\zeta\mathbf{G}:=\sum _{k=1}^N\zeta _kG_k$ is an
analogue of the system matrix $G$ (see \eqref{eq:matr}). Moreover, the
pencil of
contractions $\zeta\mathbf{A}$ (resp. $\zeta\mathbf{G}$) is an analogue
of the contraction $A$ (resp. $G$), and the pencil of unitary operators
$\zeta\mathbf{G}$ (of unitary colligations) is an analogue of the
unitary operator $G$ (of the unitary colligation)
($\zeta\in\nspace{T}{N}$). (c$'$)~The multiparametric semigroup of
operators which reproduces the dynamics of a system and serves as an
analogue of the Lax-Phillips semigroup associated with a one-parametric
system is constructed. (d$'$)~The notion of dilation for
multiparametric
systems is introduced, and by means of it some multidimensional
analogues of results in control theory as well as the criterion of the
existense of conservative dilations for multiparametric dissipative
scattering LSDSs (in particular, of unitary dilations for linear pencils
of contractive operators in a Hilbert space that are considered
on $\nspace{T}{N}$) are obtained.

Let us notice that generally speaking there is an extensive bibliography
on multiparametric or, as one usually says, multidimensional systems
theory (see e.g. a survey \cite{Bo}). We shall indicate here only those
approaches where metric properties of systems (dissipativity, conservativity
and their analogues) are essential and where some of aspects extracted
above of systems theory are reflected. First of all, it is a system
approach to the investigation of $N$-tuples of commuting non-selfadjoint
operators in a Hilbert space having its origin in works of
M.S.~Liv\v{s}ic and his associates (see \cite{LI}, \cite{LW}, and also
\cite{LKMV} and bibliography indicated there). In this approach systems
with continuous ``multidimensional time'' are overdetermined that brings
to the necessity of additional relations which have to be imposed upon
the $N$-tuple of the main operators of such systems. As a
consequence this induces connections with the theory of functions on
Riemann surfaces (in fact, not with functions of several complex
variables). In the discrete case, as V. Vinnikov told to the author, in
his joint work with J.A.~Ball an analogous systems theory is
constructed, that is connected with the investigation of $N$-tuples of
commuting contractions in a Hilbert space and also brings to functions
on Riemann surfaces (these results are unpublished yet).
Multidimensional analogues of dissipative resistance systems (see \cite{A2})
and their connections with the theory of functions of several complex
variables were considered in \cite{Be}. The different multidimensional
generalizations of the Lax-Phillips scattering scheme in the continuous
case were constructed in \cite{Se} and \cite{Z}. In the discrete case
the multidimensional analogue of the abstract scattering scheme
(including not only the Lax-Phillips scheme \cite{LP} but also the
Adamjan-Arov scheme \cite{AA}) has appeared in \cite{CS}. The structure
of a multiparametric semigroup of isometries was investigated in
\cite{Su} and \cite{GS}. Finally, in the paper by J.A.~Ball and
T.T.~Trent \cite{BT} systems known in literature on the multidimensional
systems theory as the ``Roesser model'' (see e.g. \cite{K}), with imposed
metric constraints, namely
conservative (in particular, unitary) systems are investigated (in the
terminology of the authors such a system is called conservative if it
satisfies one of ``energy'' equalities, and unitary if it satisfies
both of these equalities analogous to \eqref{eq:unitary}). Using the
result of
\cite{Ag} on the realization of operator-valued functions from
the generalized Schur class (we shall recall its definition in
Section~\ref{sec:class}) by unitary systems of the above-mentioned type,
in
\cite{BT}
various functional models of those systems are constructed, and also
interpolation problems for holomorphic functions on $\nspace{D}{N}$ are
solved.

Our paper has the following structure. In Section~\ref{sec:n-sys} the
appropriate for
the further generalization reformulation of the definition of LSDS with
(one-dimensional) discrete time and corresponding reformulations of main
results on such systems are given. Then the definition of a
multiparametric LSDS and related definitions of the conjugate system,
the transfer function, and also of the associated Lax-Phillips semigroup
and of the associated one-parametric LSDS that reproduce the
dynamics of such a system are given.
In Section~\ref{sec:dissipative-n-sys}
multiparametric dissipative, and conservative scattering LSDSs are defined
and their characterizations in terms of the associated
Lax-Phillips semigroup and of the associated one-parametric system are
obtained. In Section~\ref{sec:class} the class of transfer functions of
multiparametric conservative scattering LSDSs is described as the
subclass of functions from the generalized Schur class that are equal to
zero at $z=0$ . Then the refinement (in one way part) of this result is
obtained: the theorem on the closely connected conservative realization
of such a function is proved.

In view of restrictions on the size of a paper we plan to present the
definition of a dilation for multiparametric systems and results on
conservative dilations of multiparametric dissipative scattering LSDSs in
our next paper.

\section{Multiparametric LSDS and its transfer
function}\label{sec:n-sys}
\subsection{Some remarks on LSDSs with (one-dimensional) discrete
time}\label{sec:rem-sys}
For the further generalization of the notion of LSDS we need to make
some renotations in \eqref{eq:sys}. Namely, for all
$t\in\mathbb{Z_+}:=\{
t\in\mathbb{Z}:t\geq0\} $ we set $\phi ^-(t):=u(t),\quad \phi
^+(t+1):=y(t)$. Then one can rewrite \eqref{eq:sys} as follows:
\begin{equation}\label{eq:sys'}
\alpha ^0:\left\{ \begin{array}{lll}
x(t)&=&Ax(t-1)+B\phi ^-(t-1),\\
\phi ^+(t)&=&Cx(t-1)+D\phi ^-(t-1),
\end{array}
\right.\quad
(t\in\mathbb{Z_+}\setminus\{ 0\} )
\end{equation}
and the initial condition \eqref{eq:init-sys} remains. We call $x(t)\
(\in\Hspace{X}),\quad \phi ^-(t)\ (\in\Hspace{N^-}),\quad \phi ^+(t+1)\
(\in\Hspace{N^+})$ for $t\in\mathbb{Z_+}$ \emph{states, input data and
output data} of the system \eqref{eq:sys'}\&\eqref{eq:init-sys},
respectively. For $\alpha ^0$
of the form \eqref{eq:sys'} we shall use also the short notation $\alpha
^0 =(A,
B, C, D; \Hspace{X}, \Hspace{N^-}, \Hspace{N^+})$. The dissipativity
condition for $\alpha ^0$ can be rewritten in the following way:
\begin{eqnarray}
\forall x(0)\in\Hspace{X},\quad \forall t\in\mathbb{Z_+}\setminus \{ 0\}
,\quad
\forall\{ \phi ^-(\tau)\mid 0\leq\tau <t\} \subset\Hspace{N^-}\nonumber
\\ {\|\phi ^-(t-1)\|}^2-{\|\phi ^+(t)\|}^2\geq{\| x(t)\|}^2-{\| x(t-
1)\|}^2,\label{eq:energ'}
\end{eqnarray}
and the conservativity conditions for $\alpha ^0$ are obtained if we
change the symbol ``$\geq$'' in \eqref{eq:energ'} by ``='' and require
the
analogous
equalities for the \emph{conjugate system} $(\alpha ^0)^*:=(A^*, C^*,
B^*,
D^*; \Hspace{X}, \Hspace{N^+}, \Hspace{N^-})$. If one interprets input
data $\phi ^-(t)$ and output data $\phi ^+(t)$ as the ``amplitudes of
the
incident and reflected waves'', ${\| \phi ^-(t)\| }^2$ and ${\| \phi
^+(t)\| }^2$ as their ``powers'', and ${\| x(t)\| }^2$ as the ``energy''
of
inner states $x(t)$, then the condition \eqref{eq:energ'} 
for $\alpha ^0$ means
the dissipation of energy, and the conservativity conditions for 
$\alpha
^0$ means the full balance of energy, i.e. its conservation both for
$\alpha ^0$ and for $(\alpha ^0)^*$. The last one may be interpreted as
a
system ``with inverse time and the inverse direction of waves
propagation''.

Substitute by virtue of a system the expression $Ax(t-2)+B\phi ^-
(t-2)$ instead of $x(t-1)$ in the right-hand side of \eqref{eq:sys'},
after this
$Ax(t-3)+B\phi ^-(t-3)$
instead of $x(t-2)$, etc. Then we get for $t\in\mathbb{Z_+}\setminus \{
0\}$
\begin{eqnarray}
x(t) & = & A^tx(o)+\sum _{\tau =0}^{t-1}A^{t-\tau -1}B\phi ^-(\tau
), \label{eq:state'}\\
\phi ^+(t) & = & CA^{t-1}x(0)+D\phi ^-(t-1)+\sum _{\tau =0}^{t-2}CA^{t-
\tau  -2}B\phi ^-(\tau)\label{eq:output'}
\end{eqnarray}
(for $t=1$ only two summands remain in \eqref{eq:output'}). If we
specify the zero
initial condition in \eqref{eq:init-sys} (i.e., put $x_0=0$) then
\begin{eqnarray}
x(t) & = & \sum _{\tau =0}^{t-1}A^{t-\tau -1}B\phi ^-(\tau ), \quad
(t\in\mathbb{Z_+}\setminus \{ 0\} ) \label{eq:specstate'}\\
\phi ^+(t) & = & D\phi ^-(t-1)+\sum _{\tau =0}^{t-2}CA^{t-\tau  -2}B\phi
^-(\tau). \quad   (t\in\mathbb{Z_+}\setminus \{ 0\}
)\label{eq:specoutput'}
\end{eqnarray}
Consider the formal power series
\begin{eqnarray}\label{eq:series-sys'}
\widehat{x}(z)=\sum _{t\in\mathbb{Z_+}\setminus \{ 0\} }x(t)z^t,\ \
\widehat{\phi}^-(z)=\sum _{t\in\mathbb{Z_+} }\phi ^-(t)z^t,\ \
\widehat{\phi}^+(z)=\sum _{t\in\mathbb{Z_+}\setminus \{ 0\}
}\phi ^+(t)z^t.
\end{eqnarray} Then it is easy to obtain (formally) from
\eqref{eq:specstate'} and \eqref{eq:specoutput'}
\begin{eqnarray}\label{eq:transf-state'}
\widehat{x}(z)&=&\left(\sum _{n=0}^\infty
{(zA)}^nzB\right)\widehat{\phi}^-(z),\\
\widehat{\phi}^+(z)&=&\left(zD+\sum _{n=0}^\infty
zC{(zA)}^nzB\right)\widehat{\phi}^-(z).\label{eq:transf-output'}
\end{eqnarray}
Since the operator $A$ is bounded, for $z\in\mathbb{C}$ from a small
neighbourhood of zero the series in \eqref{eq:transf-state'} and
\eqref{eq:transf-output'} are convergent in
the operator norm (moreover, this covergence is uniform on compact sets
in this neighbourhood). Thus the operator-valued functions
\begin{eqnarray*}
{(I_\Hspace{X}-zA)}^{-1}zB&=&\sum _{n=0}^\infty {(zA)}^nzB,\\
zD+zC{(I_\Hspace{X}-zA)}^{-1}zB&=&zD+\sum _{n=0}^\infty zC{(zA)}^nzB
\end{eqnarray*}
turn out to be holomorphic on this neighbourhood. We shall call the
operator-valued function
\begin{equation}\label{eq:tf-sys'}
\theta _{{\alpha }^0}(z):=zD+zC{(I_\Hspace{X}-zA)}^{-1}zB
\end{equation}
the \emph{transfer function} of a system $\alpha ^0$ of the form
\eqref{eq:sys'}.

>From \eqref{eq:transf-state'} and \eqref{eq:transf-output'} one can
deduce the formal relations
\begin{equation}\label{eq:transform-sys'}
\widehat{\alpha ^0}:\left\{ \begin{array}{lll}
\widehat{x}(z)&=&zA\widehat{x}(z)+zB\widehat{\phi}^-(z),\\
\widehat{\phi}^+(z)&=&zC\widehat{x}(z)+zD\widehat{\phi}^-(z).
\end{array}
\right.
\end{equation}
As in \cite{BC} for systems of the form \eqref{eq:sys}, we shall call
$\widehat{\alpha ^0}$ the $Z$-\emph{transform} of $\alpha ^0$. If the
$\Hspace{N^-}$-valued function $\widehat{\phi}^-(z)$ from
\eqref{eq:series-sys'} is
holomorphic on some neighbourhood of zero in $\mathbb{C}$ then by
\eqref{eq:transf-state'} and \eqref{eq:transf-output'}
$\widehat{x}(z)$ and
$\widehat{\phi}^+(z)$ from \eqref{eq:series-sys'} are also holomorphic
on a neighbourhood
of $z=0\quad\Hspace{X}$-valued (resp. $\Hspace{N^+}$-valued) functions,
thus \eqref{eq:transform-sys'} turns out to be a system of equations
with holomorphic
functions, \eqref{eq:transf-state'} and \eqref{eq:transf-output'} turn
into equations
\begin{eqnarray}\label{eq:transf-state'1}
\widehat{x}(z)&=&{(I_\Hspace{X}-zA)}^{-1}zB\widehat{\phi}^-(z),\\
\widehat{\phi}^+(z)&=&\theta _{\alpha ^0}(z)\widehat{\phi}^-
(z),\label{eq:transf-output'1}
\end{eqnarray}
that are equivalent to (1.11) in this case.

Comparing \eqref{eq:tf-sys'} with \eqref{eq:tf-sys} we obtain
\begin{equation}\label{eq:tf-tf'}
\theta _{\alpha ^0}(z)=z\theta _\alpha (z).
\end{equation}
Let us define the class $S^0(\Hspace{N^-}, \Hspace{N^+})$ as the
subclass of all operator-valued functions from $S(\Hspace{N^-},
\Hspace{N^+})$ that are equal to zero at $z=0$. If $\theta\in
S(\Hspace{N^-}, \Hspace{N^+})$ then $z\theta\in S^0(\Hspace{N^-},
\Hspace{N^+})$. Conversely, if $\psi\in S^0(\Hspace{N^-}, \Hspace{N^+})$
then by the Schwarz lemma for Banach space-valued functions (see e.g.
Section 8.1.2 in \cite{R}) $\psi =z\theta$ with $\theta\in S(\Hspace{N^-
}, \Hspace{N^+})$. So we have the canonical bijection between the
classes $S(\Hspace{N^-}, \Hspace{N^+})$ and $S^0(\Hspace{N^-},
\Hspace{N^+})$:
\begin{equation}\label{eq:schwarz}
S^0(\Hspace{N^-}, \Hspace{N^+})=zS(\Hspace{N^-}, \Hspace{N^+}).
\end{equation}
From
\eqref{eq:tf-tf'} and \eqref{eq:schwarz} we obtain the following
reformulations of results
cited in the previous section.
\begin{thm}\label{thm:tf-sys}
The transfer function $\theta_{\alpha ^0}$
of an arbitrary dissipative scattering LSDS $\alpha ^0$ of the form
\eqref{eq:sys'}
belongs to the class $S^0(\Hspace{N^-}, \Hspace{N^+})$.
\end{thm}
\begin{thm}\label{thm:cc-cons}
Any function $\theta\in S^0(\Hspace{N^-}, \Hspace{N^+})$
allows a closely connected conservative realization of the form
\eqref{eq:sys'}.
This realization is unique up to unitary similarity.
\end{thm}
\begin{rem}\label{rem:similar}
The definition of unitary similarity for systems \eqref{eq:sys'} is the
same as
for systems \eqref{eq:sys} (see e.g. \cite{BC}).
\end{rem}
Theorem~\ref{thm:tf-sys} and Theorem~\ref{thm:cc-cons} imply the
following.
\begin{thm}\label{thm:cons-realiz}
The class of transfer functions of
conservative
scattering LSDSs with the input space $\Hspace{N^-}$ and the output
space $\Hspace{N^+}$ coincides with $S^0(\Hspace{N^-}, \Hspace{N^+})$.
\end{thm}
One can correspond to a system $\alpha ^0$ of the form \eqref{eq:sys'}
(as well as
to a system $\alpha $ of the form \eqref{eq:sys} the \emph{associated
Lax-Phillips semigroup} (see \cite{AA}, \cite{A1},
\cite{N}) $\mathfrak{W}_{\alpha ^0}=\{ W_{{\alpha
}^0}^t\mid t\in\mathbb{Z_+}\}$ where
$W_{\alpha ^0}\in [\Hspace{H}_{{\alpha }^0}, \Hspace{H}_{{\alpha }^0}]$
is its \emph{generator}, $\Hspace{H}_{{\alpha
}^0}:=\Hspace{D_+}\oplus\Hspace{X}\oplus\Hspace{D_-
},\quad\Hspace{D_+}:=\cdots
\oplus\Hspace{N^+}\oplus\Hspace{N^+},\quad\Hspace{D_-}:=\Hspace{N^-}
\oplus\Hspace{N^-}\oplus\cdots .$ If $h\in \Hspace{H}_{{\alpha }^0}$,
i.e.,\begin{equation}\label{eq:col-vect}
h=\mbox{col}(\ldots, v_{-1}, v_0; \fbox{\ensuremath{x_0}}; u_0,
u_1,\ldots )
\end{equation}
(here the element $x_0$ of the subspace $\Hspace{X}$ in
$\Hspace{H}_{\alpha ^0}$ is distinguished by a frame), then
\begin{displaymath}
W_{\alpha ^0}h:=\mbox{col}(\ldots , v_0, Cx_0+Du_0;
\fbox{\ensuremath{Ax_0+Bu_0}};u_1,
u_2,\ldots ).
\end{displaymath}
Evidently, this semigroup fits in the abstract scattering scheme by Lax-
Phillips which was described in Introduction. For the conjugate system
$(\alpha ^0)^*$ the associated semigroup
$\mathfrak{W}_{(\alpha ^0)^*}$
is a ``conjugate semigroup with inverse time'', i.e., for
$h\in\Hspace{H}_{\alpha ^0}$ from \eqref{eq:col-vect} one defines the
transform
\begin{displaymath}
\gamma h:=\mbox{col}(\ldots , u_1, u_0; \fbox{\ensuremath{x_0}};
v_0, v_{-1},
\ldots )
\end{displaymath}
mapping the space $\Hspace{H}_{\alpha ^0}$ isometrically onto
$\Hspace{H}_{(\alpha ^0)^*}$, and then
\begin{displaymath}
W_{(\alpha ^0)^*}^t=\gamma {(W_{\alpha ^0}^*)}^t\gamma ^{-1}.\quad (t\in
\mathbb{Z_+})
\end{displaymath}

The semigroup $\mathfrak{W}_{\alpha ^0}$ is contractive (resp. unitary)
if and only if $\alpha ^0$ (as well as $\alpha $) is a dissipative (resp.
conservative) scattering LSDS. The semigroup $\mathfrak{W}_{\alpha ^0}$
reproduces the dynamics of the system $\alpha ^0$ in the following
sense: if the sequence $\{ u_\tau\mid\tau =0,1,\ldots \}$ from
\eqref{eq:col-vect} is
supplied to the input of $\alpha ^0$, i.e., $\phi ^-(\tau )=u_\tau
\quad (\tau =0,1,\ldots )$,
and $x_0$ from \eqref{eq:col-vect}
is
substituted into the initial condition \eqref{eq:init-sys} then for $h$
from \eqref{eq:col-vect} we
have:
\begin{displaymath}
W_{\alpha ^0}^th=\mbox{col}(\ldots , v_{-1}, v_0, \phi ^+(1),
\phi ^+(2),
\ldots ,
\phi ^+(t); \fbox{\ensuremath{x(t)}}; \phi ^-(t),\phi ^-(t+1),\ldots )
\end{displaymath}
where $\phi ^+(1),\ \phi ^+(2), \ldots , \phi ^+(t),\ x(t)$
turn out to be the output data and states of the system $\alpha ^0$ at
corresponding moments of time $t\in \mathbb{Z_+}\setminus \{ 0\}$.

\subsection{The definition of a multiparametric LSDS}\label{sec:def-n-
sys}

Let us introduce the notion of multiparametric LSDS generalizing the
notion of LSDS \eqref{eq:sys'}\&\eqref{eq:init-sys}. For $t\in
\nspace{Z}{N}$ we set $|t|:=\sum
_{k=1}^Nt_k,\quad \nspace{\widetilde{Z}}{N}_+:=\{ t\in
\nspace{Z}{N}:|t|\geq 0\} ,\quad \nspace{\widetilde{Z}}{N}_0:=\{ t\in
\nspace{Z}{N}:|t|=0\}$, for $k\in\{1,\ldots ,N\} \quad e_k:=(0,\ldots ,
0, 1, 0,\ldots ,0)\in \nspace{Z}{N}$
(here unit is on the $k$-th place and zeros are otherwise). We define a
\emph{multiparametric LSDS} as the following system of equalities:
\begin{equation}\label{eq:nsys}
\alpha :\left\{\begin{array}{lll}
x(t)&=&\sum _{k=1}^N(A_kx(t-e_k)+B_k\phi ^-(t-e_k)),\\
\phi ^+(t)&=&\sum _{k=1}^N(C_kx(t-e_k)+D_k\phi ^-(t-e_k)),
\end{array}
\right. (t\in\nspace{\widetilde{Z}}{N}_+\setminus
\nspace{\widetilde{Z}}{N}_0)
\end{equation}
where for all $t\in\nspace{\widetilde{Z}}{N}_+\quad x(t)\
(\in\Hspace{X}),\quad \phi ^-(t)\ (\in\Hspace{N^-})$, and for all
$t\in\nspace{\widetilde{Z}}{N}_+\setminus
\nspace{\widetilde{Z}}{N}_0\quad \phi ^+(t)\
(\in\Hspace{N^+})$ are respectively \emph{states},  \emph{input data}
and  \emph{output data}
of $\alpha $, and $\Hspace{X},\ \Hspace{N^-},\ \Hspace{N^+}$ are
separable Hilbert spaces that are called the \emph{state space}, the
\emph{input space} and the \emph{output space} respectively, for all
$k\in\{1,\ldots ,N\} \ A_k\in [\Hspace{X}, \Hspace{X}],\quad B_k\in
[\Hspace{N^-},\Hspace{X}],\quad C_k\in [\Hspace{X},\Hspace{N^+}],\quad
D_k\in [\Hspace{N^-},\Hspace{N^+}]$, and also the following analogue of
the initial condition \eqref{eq:init-sys} is given:
\begin{equation}\label{eq:init-nsys}
x(t)=x_0(t),\quad (t\in\nspace{\widetilde{Z}}{N}_0)
\end{equation}
where $x_0(\cdot ):\nspace{\widetilde{Z}}{N}_0\to\Hspace{X}$ is a
prescribed function. If one denotes the $N$-tuple of operators $T_k\quad
(k=1,\ldots ,N)$ by $\mathbf{T}:=(T_1,\ldots ,T_N)$ then for such a
system one may use the short notation $\alpha =(N; \mathbf{A},
\mathbf{B}, \mathbf{C}, \mathbf{D};\Hspace{X}, \Hspace{N^-},
\Hspace{N^+})$. It is clear that
$(1; \mathbf{A}, \mathbf{B},
\mathbf{C},\mathbf{D};\Hspace{X},\Hspace{N^-}, \Hspace{N^+})=(A, B, C,
D; \Hspace{X}, \Hspace{N^-}, \Hspace{N^+})$.
\begin{rem}\label{rem:octant}
A system $\alpha $ can be considered on the positive octant
$\nspace{Z}{N}_+:=\{ t\in\nspace{Z}{N}:t_k\geq 0,\quad k=1,\ldots
,N\}$ only, i.e., we can take an input function $\phi ^-(\cdot)$ with
the support in $\nspace{Z}{N}_+$ and choose the following initial data
in \eqref{eq:init-nsys}:
\begin{displaymath}
x_0(t)=\left\{\begin{array}{ll}
0    & \mbox{for $t\in\nspace{\widetilde{Z}}{N}_0\setminus \{ 0\}$},\\
x_0  & \mbox{for $t=0$},
\end{array}
\right.
\end{displaymath}
with some prescribed $x_0\in \Hspace{X}$, and then according to
\eqref{eq:nsys} \
$\mbox{supp\ }x(\cdot)\subset\nspace{Z}{N}_+,\quad
\mbox{supp\ } \phi
^+(\cdot)\subset\nspace{Z}{N}_+$.

\end{rem}

Substitute by virtue of a system the expression $\sum _{j=1}^N(A_jx(t-
e_k-e_j)+B_j\phi ^-
(t-e_k-e_j))$ instead of $x(t-e_k)$ in the right-hand side of
\eqref{eq:nsys},
after this the corresponding expression
instead of $x(t-e_k-e_j)$, etc. Then we get for
$t\in\nspace{\widetilde{Z}}{N}_+\setminus
\nspace{\widetilde{Z}}{N}_0$
\begin{eqnarray*}
x(t) & = &
      \sum _{k_1=1}^N\cdots\sum _{k_{|t|}=1}^NA_{k_1}\cdots
A_{k_{|t|}}x(t-\sum _{j=1}^{|t|}e_{k_j})\\
     & \mbox{} +& \sum _{l=1}^{|t|}\sum _{k_{1}=1}^N\cdots\sum
_{k_{l-1}=1}^N \sum
_{k_l=1}^NA_{k_1}\cdots A_{k_{l-1}}B_{k_l}\phi ^-(t-\sum
_{j=1}^le_{k_j}),\\
\phi ^+(t) & = &
 \sum _{k_1=1}^N\sum _{k_2=1}^N\cdots\sum
_{k_{|t|}=1}^NC_{k_1}A_{k_2}\cdots
A_{k_{|t|}}x(t-\sum _{j=1}^{|t|}e_{k_j})+\sum _{k=1}^ND_k\phi ^-(t-
e_k)\\
& \mbox{} + & \sum _{l=2}^{|t|}\sum _{k_{1}=1}^N\sum
_{k_{2}=1}^N\cdots\sum
_{k_{l-1}=1}^N \sum
_{k_l=1}^NC_{k_1}A_{k_2}\cdots A_{k_{l-1}}B_{k_l}\phi ^-(t-\sum
_{j=1}^le_{k_j}).
\end{eqnarray*}
Let $\tau \leq t$ mean $t-\tau\in \nspace{Z}{N}_+$. We denote by
\begin{displaymath}
c_s:=\frac{|s|!}{s_1!\cdots s_N!}\quad
(s\in\nspace{Z}{N}_+)
\end{displaymath}
the numbers of permutations with repetitions (the polynomial
coefficients). We introduce also the following
notations:\begin{eqnarray}
\mathbf{A}^s &:=&
c_s^{-1}\sum _\sigma A_{[\sigma
(1)]}\cdots
A_{[\sigma (|s|)]},\quad (s\in \nspace{Z}{N}_+) \label{eq:multipower} \\
{(\mathbf{A\sharp
B})}^s& :=& c_s^{-1}\sum _\sigma A_{[\sigma
(1)]}\cdots A_{[\sigma (|s|-1)]}B_{[\sigma (|s|)]}, (s\in
\nspace{Z}{N}_+\setminus
\{ 0\} )\label{eq:sharp} \\
{(\mathbf{C\flat A})}^s&:=& c_s^{-1}\sum
_\sigma C_{[\sigma
(1)]}A_{[\sigma (2)]}\cdots A_{[\sigma (|s|)]},\quad (s\in
\nspace{Z}{N}_+\setminus
\{ 0\} )\label{eq:flat}\\
{(\mathbf{C\flat A\sharp B})}^s& :=&c_s^{-1}\sum _\sigma C_{[\sigma
(1)]}A_{[\sigma (2)]}\cdots A_{[\sigma
(|s|-1)]}B_{[\sigma (|s|)]}\nonumber \\
&&(s\in \nspace{Z}{N}_+\setminus \{ 0,
e_1,\ldots, e_N\} )\label{eq:flat-sharp}
\end{eqnarray}
for \emph{symmetrized multipowers of the operator} $N$-\emph{tuple}
$\mathbf{A}$~\eqref{eq:multipower}, \emph{of the } $N$-tuple
$\mathbf{A}$
\emph{bordered with the} $N$-tuple $\mathbf{B}$ \emph{from the
right}~\eqref{eq:sharp}, \emph{of the } $N$-tuple $\mathbf{A}$
\emph{bordered with
the} $N$-tuple $\mathbf{C}$ \emph{from the left}~\eqref{eq:flat},
\emph{of the }
$N$-tuple $\mathbf{A}$ \emph{bordered with the} $N$-tuple $\mathbf{C}$
\emph{from the left and with the} $N$-tuple $\mathbf{B}$ \emph{from the
right}~\eqref{eq:flat-sharp}. In these formulas the summation index
$\sigma$ runs the
set of all permutations of $|s|$ elements of $N$ different types with
repetitions (an element of the $j$-th type repeats itself $s_j$ times,
$[k]\ (\in\{1,\ldots ,N\})$ denotes the type of an element $k$). Note
that in the case of a commutative $N$-tuple $\mathbf{A}$ we have
$\mathbf{A}^s=A_1^{s_1}\cdots A_N^{s_N}$, i.e., a usual multipower. In
notations \eqref{eq:multipower}--\eqref{eq:flat-sharp} we obtain for
$t\in\nspace{\widetilde{Z}}{N}_+\setminus
\nspace{\widetilde{Z}}{N}_0$
\begin{eqnarray}
x(t)& =& \sum _{\tau\leq t,\ |\tau |=0}c_{t-\tau }\mathbf{A}^{t-\tau
}x(\tau )+\sum _{\tau\leq t,\ \tau\neq t}c_{t-\tau }{(\mathbf{A\sharp
B})}^{t-\tau }\phi ^-(\tau ),\label{eq:n-state}\\
\phi ^+(t)& =& \sum _{\tau\leq t,\ |\tau |=0}c_{t-\tau }{(\mathbf{C\flat
A})}^{t-\tau }x(\tau
)+\sum _{k=1}^ND_k\phi ^-(t-e_k)\nonumber\\
         & \mbox{} +& \sum _{\tau\leq t,\ |t-\tau |\geq
2}c_{t-\tau }{(\mathbf{C\flat A \sharp B})}^{t-\tau }\phi ^-(\tau ).
\label{eq:n-output}
\end{eqnarray}
It is easy to assure oneself that for $N=1$ \eqref{eq:n-state} and
\eqref{eq:n-output} coincide
with \eqref{eq:state'} and \eqref{eq:output'} respectively.

\subsection{The $Z$-transform and the transfer function}\label{sec:tf}
Let us specify the zero condition in \eqref{eq:init-nsys}, i.e.,
$x_0(t)=0$
for all $t\in\nspace{\widetilde{Z}}{N}_0$. Then from \eqref{eq:n-state}
and
\eqref{eq:n-output} we get for $t\in\nspace{\widetilde{Z}}{N}_+\setminus
\nspace{\widetilde{Z}}{N}_0$
\begin{eqnarray}
x(t)      & = &\sum _{\tau\leq t,\ \tau\neq t}c_{t-\tau
}{(\mathbf{A\sharp
 B})}^{t-\tau }\phi ^-(\tau ), \label{eq:n-specstate}\\
\phi ^+(t)& = &\sum _{k=1}^ND_k\phi ^-(t-e_k)+\sum _{\tau\leq t,\ |t-
\tau
|\geq 2}c_{t-\tau }{(\mathbf{C\flat A \sharp B})}^{t-\tau }\phi ^-(\tau
).\label{eq:n-specoutput}
\end{eqnarray}
Consider the formal power series
\begin{eqnarray}
\widehat{x}(z)=\sum _{t\in\nspace{\widetilde{Z}}{N}_+\setminus
\nspace{\widetilde{Z}}{N}_0}x(t)z^t,\quad
\widehat{\phi}^-(z)=\sum _{t\in\nspace{\widetilde{Z}}{N}_+}\phi ^-
(t)z^t,\nonumber\\
\widehat{\phi}^+(z)=\sum _{t\in\nspace{\widetilde{Z}}{N}_+\setminus
\nspace{\widetilde{Z}}{N}_0}\phi ^+(t)z^t \label{eq:n-series}
\end{eqnarray}
(here $z^t:=z_1^{t_1}\cdots z_N^{t_N}$ for any $N$-tuple $z=(z_1,\ldots
, z_N)$ of commuting variables and for any $t\in \nspace{Z}{N}_+$).
Then from \eqref{eq:n-specstate} we get (formally)
\begin{eqnarray*}
\lefteqn{\widehat{x}(z) =\sum _{t\in\nspace{\widetilde{Z}}{N}_+\setminus
\nspace{\widetilde{Z}}{N}_0}z^t\sum _{\tau\leq t,\ \tau\neq t}c_{t-\tau
}{(\mathbf{A\sharp
 B})}^{t-\tau }\phi ^-(\tau )}\\
               & &=\sum _{s\in \nspace{Z}{N}_+\setminus \{ 0\} }
c_sz^s{(\mathbf{A\sharp
 B})}^s\sum _{\tau\in\nspace{\widetilde{Z}}{N}_+}z^\tau\phi ^-(\tau )=
\sum _{n=0}^\infty {\left(\sum _{k=1}^Nz_kA_k\right)}^n\sum
_{j=1}^Nz_jB_j\widehat{\phi }^-(z).
\end{eqnarray*}
Using the notation $z\mathbf{T}:=\sum _{k=1}^Nz_kT_k$
for $N$-tuples $z=(z_1,\ldots
, z_N)$ and $\mathbf{T}=(T_1, \ldots
,T_N) $ we get
\begin{equation}\label{eq:transf-n-state}
\widehat{x}(z)=\left(\sum _{n=0}^\infty
{(z\mathbf{A})}^nz\mathbf{B}\right)\widehat{\phi }^-(z)
\end{equation}
(cf. \eqref{eq:transf-state'}). From \eqref{eq:n-specoutput} we get
(formally)
\begin{eqnarray*}
\lefteqn{\widehat{\phi }^+(z) = \sum
_{t\in\nspace{\widetilde{Z}}{N}_+\setminus
\nspace{\widetilde{Z}}{N}_0}z^t\left(\sum _{k=1}^ND_k\phi ^-(t-e_k)+\sum
_{\tau\leq t,\ |t-\tau
|\geq 2}c_{t-\tau }{(\mathbf{C\flat A \sharp B})}^{t-\tau }\phi ^-(\tau
)\right)}\\
                     &&= \sum _{k=1}^Nz_kD_k\sum
_{\tau\in\nspace{\widetilde{Z}}{N}_+}z^\tau\phi ^-(\tau )+\sum
_{s\in\nspace{Z}{N}_+\setminus \{ 0, e_1, \ldots ,e_N\} }c_sz^s
{(\mathbf{C\flat A \sharp B})}^s\sum
_{\tau\in\nspace{\widetilde{Z}}{N}_+}z^\tau\phi ^-(\tau ).
\end{eqnarray*}
Finally,
\begin{equation}\label{eq:transf-n-output}
\widehat{\phi }^+(z)=\left(z\mathbf{D}+\sum _{n=0}^\infty
z\mathbf{C}{(z\mathbf{A})}^nz\mathbf{B}\right)\widehat{\phi}^-(z)
\end{equation}
(cf. \eqref{eq:transf-output'}). Since $\mathbf{A}$ is a $N$-tuple of
bounded operators, in
a small neighbourhood of $z=0$ in $\nspace{C}{N}$ the series in
\eqref{eq:transf-n-state}
and \eqref{eq:transf-n-output} are convergent in the operator norm
(moreover, this
convergence is uniform on compact sets in this neighbourhood). Thus the
operator-valued functions
\begin{eqnarray*}
{(I_\Hspace{X}-z\mathbf{A})}^{-1}z\mathbf{B}&=&\sum _{n=0}^\infty
{(z\mathbf{A})}^nz\mathbf{B},\\
z\mathbf{D}+z\mathbf{C}{(I_\Hspace{X}-z\mathbf{A})}^
{-1}z\mathbf{B}&=&z\mathbf{D}+\sum
_{n=0}^\infty z\mathbf{C}{(z\mathbf{A})}^nz\mathbf{B}
\end{eqnarray*}
turn out to be holomorphic on this neighbourhood. We shall call the
operator-valued function
\begin{equation}\label{eq:tf-n-sys}
\theta _\alpha (z)=z\mathbf{D}+z\mathbf{C}{(I_\Hspace{X}-z\mathbf{A})}^
{-1}z\mathbf{B}
\end{equation}
the \emph{transfer function} of a system $\alpha $ of the form
\eqref{eq:nsys}
(cf. \eqref{eq:tf-sys'}). From \eqref{eq:transf-n-state} and
\eqref{eq:transf-n-output} one can deduce the formal relations
\begin{equation}\label{eq:transform-n-sys}
\widehat{\alpha }:\left\{ \begin{array}{lll}
\widehat{x}(z) &=& z\mathbf{A}\widehat{x}(z)+z\mathbf{B}\widehat{\phi
}^-(z),\\
\widehat{\phi
}^+(z)&=&z\mathbf{C}\widehat{x}(z)+z\mathbf{D}\widehat{\phi
}^-(z)
\end{array}
\right.
\end{equation}
(cf. \eqref{eq:transform-sys'}). We shall call $\widehat{\alpha }$ the
$Z$-\emph{transform}
 of a system $\alpha $ of the form \eqref{eq:nsys}. If the $\Hspace{N}^-
$-valued
function $\widehat{\phi }^-(z)$ from \eqref{eq:n-series} is holomorphic
on some
neighbourhood of $z=0$ in $\nspace{C}{N}$ then by
\eqref{eq:transf-n-state} and \eqref{eq:transf-n-output}
$\widehat{x}(z)$ and $\widehat{\phi }^+(z)$ from
\eqref{eq:n-series} are
also holomorphic on a neighbourhood of $z=0\quad \Hspace{X}$-valued
(resp. $\Hspace{N}^+$-valued) functions, thus
\eqref{eq:transform-n-sys} turns out to be a
system of equations with holomorphic functions,
\eqref{eq:transf-n-state} and \eqref{eq:transf-n-output} turn
into equations
\begin{eqnarray}
\widehat{x}(z)&=&{(I_\Hspace{X}-z\mathbf{A})}^{-
1}z\mathbf{B}\widehat{\phi }^-(z),\label{eq:transf-n-state1}\\
\widehat{\phi }^+(z)&=&\theta _\alpha (z)\widehat{\phi }^-
(z)\label{eq:transf-n-output1}
\end{eqnarray}
(cf. \eqref{eq:transf-state'1} and \eqref{eq:transf-output'1}), that are
equivalent to \eqref{eq:transform-n-sys} in this case.

\subsection{The conjugate system}\label{sec:conj}
Let $\alpha =(N; \mathbf{A}, \mathbf{B}, \mathbf{C},
\mathbf{D};\Hspace{X}, \Hspace{N^-}, \Hspace{N^+})$ be a multiparametric
LSDS. Then we call  $\alpha ^*=(N; \mathbf{A}^*, \mathbf{C}^*,
\mathbf{B}^*, \mathbf{D}^*;\Hspace{X}, \Hspace{N^+}, \Hspace{N^-})$ the
\emph{conjugate LSDS} for $\alpha $  where for a $N$-tuple
$\mathbf{T}=(T_1, \ldots ,T_N)$ of operators from $[\Hspace{Y},
\Hspace{V}]$ we set $\mathbf{T}^*:=(T_1^*, \ldots ,T_N^*)$, that is a
$N$-tuple of operators from $[\Hspace{V}, \Hspace{Y}]$. For an arbitrary
function $\theta (z)=\theta (z_1, \ldots ,z_N)$ with values in
 $[\Hspace{Y}, \Hspace{V}]$ we define the function $\theta ^*
(z):={\theta (\bar{z})}^*={\theta (\bar{z_1}, \ldots ,\bar{z_N})}^*$
with values in $[\Hspace{V}, \Hspace{Y}]$.
\begin{prop}\label{prop:tf-conj}
Let $\alpha =(N; \mathbf{A}, \mathbf{B}, \mathbf{C},
\mathbf{D};\Hspace{X}, \Hspace{N^-}, \Hspace{N^+})$ be a multiparametric
LSDS, $\theta _\alpha (z)$ be its transfer function which is defined and
holomorphic on some open neighbourhood $\Upsilon $ of $z=0$ in
$\nspace{C}{N}$. Then the transfer function  $\theta _{\alpha ^*}(z)$ of
the conjugate LSDS $\alpha ^*$ is defined and holomorphic on $\Upsilon
^*:=\{ z\in\nspace{C}{N}:\bar{z}\in\Upsilon\} $ and $\theta _{\alpha
^*}(z)=\theta ^*_\alpha (z)$ for all $z\in \Upsilon ^*$.
\end{prop}
\begin{proof}
Evidently, if $\Upsilon$ is an open neighbourhood of $z=0$ then
$\Upsilon ^*$ is also an open neighbourhood of $z=0$. For $z\in \Upsilon
^*$ we have $\bar{z}\in\Upsilon ,\quad \theta ^*_\alpha (z)={\theta
_\alpha (\bar{z})}^*$ is holomorphic on $\Upsilon ^*$ and
\begin{eqnarray*}
\theta ^*_\alpha (z) &=& {\theta  _\alpha
(\bar{z})}^*={(\bar{z}\mathbf{D}+\bar{z}\mathbf{C}{(I_\Hspace{X}-
\bar{z}\mathbf{A})}^{-1}\bar{z}\mathbf{B})}^*\\
&= &z\mathbf{D}^*+z\mathbf{B}^*{(I_\Hspace{X}-z\mathbf{A}^*)}^{-
1}z\mathbf{C}^*=\theta _{\alpha ^*}(z)
\end{eqnarray*}
that completes the proof.
\end{proof}

\subsection{The associated semigroup and the associated one-parametric
LSDS}\label{sec:semigroup}
Let us introduce for LSDS $\alpha =(N; \mathbf{A}, \mathbf{B},
\mathbf{C}, \mathbf{D};\Hspace{X}, \Hspace{N^-}, \Hspace{N^+})$ the
multiparametric analogue of the associated Lax-Phillips semigroup (see
Section~\ref{sec:rem-sys}). Set $\nspace{\widetilde{Z}}{N}_-:=\{
t\in\nspace{Z}{N}:|t|\leq 0\}$ and let $\widetilde{\Hspace{D}}_\pm
:=l^2(\nspace{\widetilde{Z}}{N}_\mp ,\Hspace{N}^\pm ),\quad
\widetilde{\Hspace{X}}:=l^2(\nspace{\widetilde{Z}}{N}_0,\Hspace{X})$ be
 Hilbert spaces of multisequences $\{ u^\pm (t)\mid
t\in\nspace{\widetilde{Z}}{N}_\mp\}\subset \Hspace{N}^\pm ,\quad \{
y(t)\mid t\in\nspace{\widetilde{Z}}{N}_0\}\subset\Hspace{X}$
respectively, such that
\begin{displaymath}
\sum _{t\in\nspace{\widetilde{Z}}{N}_\mp}{\| u^\pm (t)\| }^2<\infty
,\quad \sum _{t\in\nspace{\widetilde{Z}}{N}_0}{\| y(t)\| }^2<\infty .
\end{displaymath}
Set also $\Hspace{H}_\alpha
:=\Hspace{\widetilde{D}}_+\oplus\Hspace{\widetilde{X}}\oplus\Hspace
{\widetilde{D}}_-$. Then we shall define the operators $W_{\alpha ,k}\in
[\Hspace{H}_\alpha ,\Hspace{H}_\alpha ]\quad (k=1,\ldots ,N)$ in the
following way. For $h=(u^+,y,u^-
)\in\Hspace{\widetilde{D}}_+\oplus\Hspace{\widetilde{X}}\oplus\Hspace{
\widetilde{D}}_-=\Hspace{H}_\alpha $ set
\begin{equation}
h_k:=W_{\alpha ,k}h=(u_k^+, y_k,u_k^-)\quad (k=1,\ldots
,N)\label{eq:gen-n-sys}
\end{equation}
where
\begin{equation}\label{eq:gen-n-sys1}
u_k^+(t)=
\left\{
\begin{array}{ll}
u^+(t+e_k), &  (|t|\leq -1)\\
\sum _{j=1}^N(C_jy(t+e_k-e_j)+D_ju^-(t+e_k-e_j)), & (|t|=0)
\end{array}
\right.
\end{equation}
\begin{eqnarray}
y_k(t)   &=&
\sum _{j=1}^N(A_jy(t+e_k-e_j)+B_ju^-(t+e_k-e_j)),\ (|t|=0)
\label{eq:gen-n-sys2}\\
u_k^-(t) &=& u^-(t+e_k).\quad  (|t|\geq 0) \label{eq:gen-n-sys3}
\end{eqnarray}
It is not difficult to assure oneself that for $k\neq j\quad W_{\alpha
,k}W_{\alpha ,j}=W_{\alpha ,j}W_{\alpha ,k}$, thus the semigroup
$\mathfrak{W}_\alpha :=\{\mathbf{W}_\alpha ^t\mid
t\in\nspace{Z}{N}_+\}\subset [\Hspace{H}_\alpha ,\Hspace{H}_\alpha ]$ is
well defined (here $\mathbf{W}_\alpha =(W_{\alpha ,1},\ldots ,W_{\alpha
,N})$ and for $t\in \nspace{Z}{N}_+\quad \mathbf{W}_\alpha ^t=W_{\alpha
,1}^{t_1}\cdots W_{\alpha ,N}^{t_N} $).
\begin{rem}\label{rem:one-param}
>From \eqref{eq:gen-n-sys}--\eqref{eq:gen-n-sys3} one can see that for
$k\in\{ 1,\ldots ,N\}\
W_{\alpha ,k}|\Hspace{\widetilde{D}}_+$ and $W_{\alpha
,k}^*|\Hspace{\widetilde{D}}_-$ are forward shift operators with
wandering generating subspaces
\begin{equation}\label{eq:wandering}
\Hspace{\widetilde{N}}^+:=\Hspace{\widetilde{D}}_+\ominus W_{\alpha
,k}\Hspace{\widetilde{D}}^+=l^2(\nspace{\widetilde{Z}}{N}_0,\Hspace{N}^+
),\  \Hspace{\widetilde{N}}^-:=\Hspace{\widetilde{D}}_-\ominus W_{\alpha
,k}^*\Hspace{\widetilde{D}}^-
=l^2(\nspace{\widetilde{Z}}{N}_0,\Hspace{N}^-)
\end{equation}
that means for $\mathfrak{W}_\alpha$ to be a multiparametric
analogue of Lax-Phillips semigroup, and for its generators $W_{\alpha
,k}\ (k=1,\ldots ,N)$ to be generators of some one-parametric 
Lax-Phillips semigroups.
\end{rem}
\begin{rem}\label{rem:reproduce}
The semigroup $\mathfrak{W}_\alpha$ reproduces the dynamics of a system
$\alpha $ in the following sense: if the multisequence $\{ u^-(t)\mid
t\in\nspace{\widetilde{Z}}{N}_+\}\ (\subset \Hspace{\widetilde{D}}_-)$
is supplied to the input of $\alpha $, i.e., $\phi ^-(t)=u^-(t)\quad
 (t\in\nspace{\widetilde{Z}}{N}_+)$, and $\{ y(t)\mid
t\in\nspace{\widetilde{Z}}{N}_0\}\ (\subset\Hspace{\widetilde{X}})$ is
substituted into \eqref{eq:init-nsys} then
\eqref{eq:gen-n-sys1}--\eqref{eq:gen-n-sys3} turn into
\begin{eqnarray*}
u_k^+(t)  &=&\left\{ \begin{array}{ll}
u^+(t+e_k), & (|t|\leq -1)\\
\phi ^+(t+e_k), & (|t|=0)
\end{array}
\right. \\
y_k(t)   &=& x(t+e_k),\quad (|t|=0)\\
u_k^-(t) &=& \phi ^-(t+e_k)\quad (|t|)\geq 0)
\end{eqnarray*}
where $x(t+e_k)$ and $\phi ^+(t+e_k)$ are the states and the output
signals of the system $\alpha $ at moments $t+e_k$. Iterating these
formulas we obtain for $h_s=(u_s^+, y_s, u_s^-):=\mathbf{W}_\alpha
^sh\quad (s\in\nspace{Z}{N}_+)$
\begin{eqnarray*}
u_s^+(t)  &=&\left\{ \begin{array}{ll}
u^+(t+s), & (|t|\leq -s)\\
\phi ^+(t+s), & (-s<|t|\leq 0)
\end{array}
\right.\\
y_s(t)   &=& x(t+s),\quad (|t|=0)\\
u_s^-(t) &=& \phi ^-(t+s)\quad (|t|)\geq 0)
\end{eqnarray*}
where $ x(t+s)$ and $\phi ^+(t+s)$ are the states and the output signals
of $\alpha $ at moments $t+s$.
\end{rem}
We shall show now that for the conjugate system $\alpha ^*$ of some
system $\alpha $ the associated semigroup $\mathfrak{W}_{\alpha ^*}$ is
a ``conjugate semigroup with inverse time''. More exactly, we define the
isomorphism $\gamma :\Hspace{H}_\alpha\to\Hspace{H}_{\alpha ^*}$ in the
following way: for $h=(u^+, y, u^-)\in\Hspace{H}_\alpha $ we set
$h^*=\gamma h:=(u_*^+, y_*, u_*^-)\in\Hspace{H}_{\alpha ^*}$ where
\begin{equation}\label{eq:inverse-t}
\left\{ \begin{array}{ll}
u_*^+(t)=u^-(-t),&(t\in\nspace{\widetilde{Z}}{N}_-)\\
y_*(t)=y(-t),&(t\in\nspace{\widetilde{Z}}{N}_0)\\
u_*^-(t)=u^+(-t)&(t\in\nspace{\widetilde{Z}}{N}_+)
\end{array}
\right.
\end{equation}
and then we have the following.
\begin{prop}\label{prop:semigr-conj}
$\forall t\in\nspace{Z}{N}_+\quad \mathbf{W}_{\alpha ^*}^t=\gamma
(\mathbf{W}_\alpha ^*)^t\gamma ^{-1}$.
\end{prop}
\begin{proof}
Evidently, it is sufficient to show that
\begin{equation}\label{eq:gen-conj}
W_{\alpha ^*,k}=\gamma W_{\alpha ,k}^*\gamma ^{-1}.\quad (k\in\{
1,\ldots
,N\} )
\end{equation}
It is not difficult to assure oneself that for $h=(u^+, y, u^-
)\in\Hspace{H}_\alpha \quad W_{\alpha ,k}^*h$ are defined by the
following formulas:
\begin{eqnarray*}
(P_\Hspace{\widetilde{D}_+}W_{\alpha ,k}^*h)(t) &=&u^+(t-e_k),\quad
(|t|\leq 0)\\
(P_\Hspace{\widetilde{X}}W_{\alpha ,k}^*h)(t) &=& \sum
_{j=1}^N(A_j^*y(t-e_k+e_j)+C_j^*u^+(t-e_k+e_j)),\quad (|t|=0)
\end{eqnarray*}
\begin{displaymath}
(P_\Hspace{\widetilde{D}_-}W_{\alpha ,k}^*h)(t) =\left\{
\begin{array}{ll}
\sum _{j=1}^N(B_j^*y(t-e_k+e_j)+D_j^*u^+(t-e_k+e_j)),& (|t|=0)\\
u^-(t-e_k), & (|t|\geq 1)
\end{array}
\right.
\end{displaymath}
hence for $h=(u^+, y, u^-)\in\Hspace{H}_\alpha$ and for $h^*=\gamma
h=(u_*^+, y_*, u_*^-)\in\Hspace{H}_{\alpha ^*}$
we have $\gamma W_{\alpha ,k}^*\gamma ^{-1}h_*=\gamma W_{\alpha
,k}^*h=:(u_{*,k}^+, y_{*,k}, u_{*,k}^-)\in \Hspace{H}_{\alpha ^*}$ where
by \eqref{eq:inverse-t}
\begin{eqnarray*}
u_{*,k}^+(t)&=&
\left\{
\begin{array}{ll}
u_*^+(t+e_k), &  (|t|\leq -1)\\
\sum _{j=1}^N(B_j^*y_*(t+e_k-e_j)+D_j^*u_*^-(t+e_k-e_j)),& (|t|=0)
\end{array}
\right.\\
y_{*,k}(t)   &=&
\sum _{j=1}^N(A_j^*y_*(t+e_k-e_j)+C_j^*u_*^-(t+e_k-e_j)),\quad (|t|=0)\\
u_{*,k}^-(t) &=& u_*^-(t+e_k).\quad  (|t|\geq 0)
\end{eqnarray*}
By the definition of $\alpha ^*$ and by
\eqref{eq:gen-n-sys}--\eqref{eq:gen-n-sys3} we get
$W_{\alpha ^*,k}h_*=(u_{*,k}^+, y_{*,k}, u_{*,k}^-)$. Taking into
account an arbitrariness of $h_*$ and unitarity of the operator
$\gamma$, \eqref{eq:gen-conj} follows.
\end{proof}
Let $\alpha =(N; \mathbf{A}, \mathbf{B}, \mathbf{C},
\mathbf{D};\Hspace{X}, \Hspace{N^-}, \Hspace{N^+})$ be a multiparametric
LSDS and $\mathfrak{W}_\alpha $ be its associated  semigroup.
Consider for an arbitrary $k\in\{ 1,\ldots
,N\}$ the one-parametric semigroup $\mathfrak{W}_{\alpha ,k}:=\{
W_{\alpha ,k}^n\mid n\in\mathbb{Z}_+\}$. In accordance with
Remark~\ref{rem:one-param}
it turns out to be the associated Lax-Phillips semigroup of some one-
parametric LSDS, namely $\alpha _k:=(\widetilde{A}_k,\widetilde{B}_k,
\widetilde{C}_k, \widetilde{D}_k; \Hspace{\widetilde{X}},
\Hspace{\widetilde{N}^-}, \Hspace{\widetilde{N}^+})$ where
$\Hspace{\widetilde{N}^+}$ and $\Hspace{\widetilde{N}^-}$
are defined in \eqref{eq:wandering},
$\widetilde{A}_k:=P_\Hspace{\widetilde{X}}W_{\alpha
,k}|\Hspace{\widetilde{X}},\quad
\widetilde{B}_k:=P_\Hspace{\widetilde{X}}W_{\alpha
,k}|\Hspace{\widetilde{N}^-},\quad
\widetilde{C}_k:=P_\Hspace{\widetilde{N}^+}W_{\alpha
,k}|\Hspace{\widetilde{X}},\quad
\widetilde{D}_k:=P_\Hspace{\widetilde{N}^+}W_{\alpha
,k}|\Hspace{\widetilde{N}^-} $. Define the unitary operators
$T_{jk}:\Hspace{\widetilde{X}}\to\Hspace{\widetilde{X}}$ and
$S_{jk}:\Hspace{\widetilde{N}^-}\to\Hspace{\widetilde{N}^-}$ for all
$j\neq k$ as follows:
\begin{eqnarray}
\forall y\in\Hspace{\widetilde{X}}\quad (T_{jk}y)(t) &:=& y(t+e_k-
e_j),\quad (t\in\nspace{Z}{N}_0) \label{eq:state-shift}\\
\forall u^-\in \Hspace{\widetilde{N}^-}\quad (S_{jk}u^-)(t) &:=&
u^-(t+e_k-e_j).\quad (t\in\nspace{Z}{N}_0)\label{eq:input-shift}
\end{eqnarray}
Then according to \eqref{eq:gen-n-sys1}--\eqref{eq:gen-n-sys3} the
system
matrix for $\alpha _k$ is
\begin{equation}\label{eq:matrices}
\widetilde{G}_k=\left (
\begin{array}{ll}
\widetilde{A}_k=A_kI_\Hspace{\widetilde{X}}+\sum _{j\neq k}A_jT_{jk} &
\widetilde{B}_k=B_kI_\Hspace{\widetilde{N}^-}+\sum _{j\neq k}B_jS_{jk}\\
\widetilde{C}_k=C_kI_\Hspace{\widetilde{X}}+\sum _{j\neq k}C_jT_{jk} &
\widetilde{D}_k=D_kI_\Hspace{\widetilde{N}^-}+\sum _{j\neq k}D_jS_{jk}
\end{array}
\right )
\end{equation}
where $A_kI_\Hspace{\widetilde{X}},\quad B_kI_\Hspace{\widetilde{N}^-
},\quad C_kI_\Hspace{\widetilde{X}},\quad D_kI_\Hspace{\widetilde{N}^-}$
are ``block diagonal'' operators, i.e.,
$(A_kI_\Hspace{\widetilde{X}}y)(t)=A_ky(t)$ for
$t\in\nspace{\widetilde{Z}}{N}_0$, etc. We shall call $\alpha _k$ the
$k$-\emph{th associated one-parametric LSDS} for $\alpha $. Each system
$\alpha _k$ reproduces the dynamics of $\alpha $: if an input sequence
of $\alpha _k$ is supplied to the input of $\alpha $, i.e., $\phi ^-
(s+ne_k)=\phi _n^-(s)$ for $s\in\nspace{\widetilde{Z}}{N}_0,\ n\in
\mathbb{Z}_+$ where $\phi ^-$ and $\phi _n^-$ are the input data of
$\alpha $ and $\alpha _k$ respectively, and the initial
condition for $\alpha _k$ is substituted into \eqref{eq:init-nsys} then
the states and
the output signals of $\alpha _k$ at moments $n>0$ coincide with
collections of the states and the output signals of $\alpha $ taken for
all $t\in\nspace{Z}{N}$ such that $|t|=n$, i.e.,
\begin{equation}\label{eq:repr}
x_n(s)=x(s+ne_k),\quad
\phi _n^+(s)=\phi ^+(s+ne_k). \quad (s\in\nspace{Z}{N}_0)
\end{equation}
Indeed, from
\eqref{eq:state-shift}--\eqref{eq:matrices} we get
for $s\in\nspace{\widetilde{Z}}{N}_0$:
\begin{eqnarray*}
x_1(s)&=&
({\widetilde{A}}_kx_0+\widetilde{B}_k\phi _0^-)(s)\\
&=&\sum _{j=1}^N(A_jx_0(s+e_k-
e_j)+B_j\phi _0^-(s+e_k-e_j))=x(s+e_k),\\
\phi _1^+(s)&=&({\widetilde{C}}_kx_0+\widetilde{D}_k\phi _0^-)(s)\\
&=&\sum _{j=1}^N(C_jx_0(s+e_k-
e_j)+D_j\phi _0^-(s+e_k-e_j))=\phi ^+(s+e_k).
\end{eqnarray*}
Iterating this calculation $n$ times we obtain \eqref{eq:repr}.

\section{Multiparametric dissipative scattering LSDS}
\label{sec:dissipative-n-sys}
\subsection{The definition and some properties}\label{sec:dissipative}
We shall call $\alpha =(N; \mathbf{A}, \mathbf{B}, \mathbf{C},
\mathbf{D};\Hspace{X}, \Hspace{N^-}, \Hspace{N^+})$
a \emph{multiparametric dissipative scattering LSDS} if for any
$\zeta\in\nspace{T}{N}$
\begin{displaymath}
\zeta\mathbf{G}:=\left (\begin{array}{ll}
\zeta\mathbf{A} & \zeta\mathbf{B}\\
\zeta\mathbf{C} & \zeta\mathbf{D}
\end{array}
\right )\in [\Hspace{X}\oplus\Hspace{N^-}, \Hspace{X}\oplus\Hspace{N^+}]
\end{displaymath}
is a contractive operator. Evidently, by virtue of the maximum principle
for holomorphic operator-valued functions (see e.g. \cite{S} ) it is
equivalent to require for the operator-valued linear function
$L_\mathbf{G}(z):=z\mathbf{G}$ to be contractive on $\nspace{D}{N}$. For
$N=1$ this notion
 coincides with the notion
 of dissipative
 scattering LSDS with discrete time (see
Introduction and Section~\ref{sec:rem-sys}).
\begin{prop}\label{prop:dissipative}
Let $\alpha =(N; \mathbf{A}, \mathbf{B}, \mathbf{C},
\mathbf{D};\Hspace{X}, \Hspace{N^-}, \Hspace{N^+})$ be given,
$\mathfrak{W}_\alpha $ and $\alpha _k$ (for some $k\in\{ 1,\ldots ,N\}$)
be the associated semigroup and the associated one-parametric LSDS
respectively. Then the following statements are equivalent:

\emph{(i)} $\alpha$ is a multiparametric dissipative scattering LSDS;

\emph{(ii)} $\mathfrak{W}_\alpha $ is a semigroup of contractions in
$\Hspace{H}_\alpha $;

\emph{(iii)} $\alpha _k$ is a dissipative scattering LSDS;

\emph{(iv)} if an input multisequence $\{\phi ^-(t)\mid
t\in\nspace{\widetilde{Z}}{N}_+\}$ of $\alpha $ satisfies 
\begin{equation}\label{eq:input-mult}
\forall n\in\mathbb{Z}_+\quad\sum _{|t|=n}{\| \phi ^-(t)\| }^2<\infty
\end{equation}
and its states $\{ x_0(t)\mid t\in
\nspace{\widetilde{Z}}{N}_0\}$ from \eqref{eq:init-nsys} satisfy 
\begin{equation}\label{eq:state-mult}
\sum _{|t|=0} {\| x_0(t)\| }^2<\infty
\end{equation}
then  $\forall n\in \mathbb{Z}_+\setminus \{ 0\}\quad \sum _{|t|=n}{\|
x(t) \| }^2<\infty ,\ \sum _{|t|=n}{\| \phi ^+(t)\| }^2<\infty$, and
\begin{equation}\label{eq:energ-dis}
\sum _{|t|=n-1}{\|\phi ^-(t)\| }^2-\sum _{|t|=n}{\|\phi ^+(t)\| }^2\geq
\sum _{|t|=n}{\| x(t)\| }^2-\sum _{|t|=n-1}{\| x(t)\| }^2.
\end{equation}
\end{prop}
\begin{proof}
(iii)$\Leftrightarrow$ (iv). If $\alpha _k$ is a dissipative scattering LSDS
then according to \eqref{eq:energ'} for any
$x_0\in\Hspace{\widetilde{X}}$ in
\eqref{eq:init-nsys}, for an arbitrary $n\in\mathbb{Z}_+\setminus \{
0\}$, and for any
finite input sequence $\{\phi _m^-\mid 0\leq
m<n\}\subset\Hspace{\widetilde{N}}^-$
we have ${\|\phi _{n-1}^-\| }^2-{\|\phi _n^+\| }^2\geq {\| x_n\| }^2-{\|
x_{n-1}\| }^2$ where $x_n\in\Hspace{\widetilde{X}}$ is the state, and
$\phi _n^+\in\Hspace{\widetilde{N}}^+$ is the output signal of $\alpha
_k$ at the moment $n$. This means that
\begin{equation}\label{eq:energ-dis-ass}
\sum _{s\in\nspace{\widetilde{Z}}{N}_0}{\|\phi _{n-1}^-(s)\| }^2-\sum
_{s\in\nspace{\widetilde{Z}}{N}_0}{\|\phi _n^+(s)\| }^2\geq
\sum _{s\in\nspace{\widetilde{Z}}{N}_0}{\| x_n(s)\| }^2-\sum
_{s\in\nspace{\widetilde{Z}}{N}_0}{\| x_{n-1}(s)\| }^2.
\end{equation}
Let $\{\phi ^-(t)\mid t\in\nspace{\widetilde{Z}}{N}_+\}$ and
$\{x_0(t)\mid
t\in\nspace{\widetilde{Z}}{N}_0\}$  satisfy
\eqref{eq:input-mult} and \eqref{eq:state-mult}, the initial state of
$\alpha _k$ be equal to $x_0=\{
x_0(s)\mid s\in\nspace{\widetilde{Z}}{N}_0\}\
(\in\Hspace{\widetilde{X}}$ according to \eqref{eq:state-mult}), and
define the input
sequence of $\alpha _k$ by $\phi _n^-(s):=\phi ^-(s+ne_k)\
(s\in\nspace{\widetilde{Z}}{N}_0)$ for all $n\in\mathbb{Z}_+$. By
\eqref{eq:input-mult}
$\phi _n^-\in\Hspace{\widetilde{N}}^-$. Then
\eqref{eq:repr} is valid.
Since for $n\in\mathbb{Z}_+\quad x_n\in\Hspace{\widetilde{X}},\ \phi
_n^+\in\Hspace{\widetilde{N}}^+$, we have
\begin{eqnarray*}
\sum _{|t|=n}{\| x(t)\| }^2=\sum _{s\in\nspace{\widetilde{Z}}{N}_0}{\|
x(s+ne_k)\| }^2=\sum _{s\in\nspace{\widetilde{Z}}{N}_0}{\| x_n(s)\| }^2=
{\| x_n\| }^2<\infty ,\\
\sum _{|t|=n}{\| \phi ^+(t)\| }^2=\sum
_{s\in\nspace{\widetilde{Z}}{N}_0}{\| \phi ^+(s+ne_k)\| }^2=\sum
_{s\in\nspace{\widetilde{Z}}{N}_0}{\| \phi^+_n(s)\| }^2=
{\| \phi^+_n\| }^2<\infty ,
\end{eqnarray*}
and \eqref{eq:energ-dis-ass} implies \eqref{eq:energ-dis}. Thus,
(iii)$\Rightarrow$(iv) holds. Reversing
the argument we obtain (iv)$\Rightarrow$(iii).

(iii)$\Leftrightarrow$(i). Define for an arbitrary separable Hilbert
space $\Hspace{Y}$ the operator
$F_\Hspace{Y}:l^2(\nspace{\widetilde{Z}}{N}_0,\Hspace{Y})\to\mathcal{L}^
2(\nspace{T}{N},\Hspace{Y})$ by
$\widehat{y}(\zeta)=(F_\Hspace{Y}y)(\zeta):=\sum
_{t\in\nspace{\widetilde{Z}}{N}_0}y(t)\zeta ^t$ (a.e.
$\zeta\in\nspace{T}{N}$) where $\mathcal{L}^2(\nspace{T}{N},\Hspace{Y})$
denotes the subspace in the space $L^2(\nspace{T}{N},\Hspace{Y})$
of square integrable on $\nspace{T}{N}\quad\Hspace{Y}$-valued
functions, that is extracted by the condition of vanishing for all
Fourier coefficients with multiindices
$t\not\in\nspace{\widetilde{Z}}{N}_0$. Evidently, $F_\Hspace{Y}$ is a
unitary operator (it is a restriction of the discrete Fourier transform
mapping $l^2(\nspace{Z}{N},\Hspace{Y})$ onto
$L^2(\nspace{T}{N},\Hspace{Y})$). If $\alpha _k$ is a dissipative scattering
LSDS then the operator ${\widetilde{G}}_k$ in \eqref{eq:matrices} is
contractive.
Thus for any $\widehat{x_0}\in \mathcal{L}^2(\nspace{T}{N},\Hspace{X})$
and
$\widehat{\phi _0^-}\in \mathcal{L}^2(\nspace{T}{N},\Hspace{N^-})$ we
have for almost every $\zeta\in\nspace{T}{N}$
\begin{eqnarray*}
\lefteqn{\left (\begin{array}{l}
\widehat{x_1}(\zeta )\\
\widehat{\phi _1^+}(\zeta )
\end{array}
\right )=\left\{\left (\begin{array}{ll}
F_\Hspace{X} & O\\
          0  & F_\Hspace{N^+}
\end{array}
\right ){\widetilde{G}}_k{\left (\begin{array}{ll}
F_\Hspace{X} & O\\
          0  & F_\Hspace{N^-}
\end{array}
\right )}^{-1}\left (\begin{array}{l}
\widehat{x_0}\\
\widehat{\phi _0^-}
\end{array}
\right )\right\} (\zeta )}\\
&&=\left (\begin{array}{l}
\sum _{t\in\nspace{\widetilde{Z}}{N}_0}\zeta ^t\sum _
{j=1}^N(A_jx_0(t+e_k-e_j)+B_j\phi _0^-(t+e_k-e_j))\\
\sum _{t\in\nspace{\widetilde{Z}}{N}_0}\zeta ^t\sum _
{j=1}^N(C_jx_0(t+e_k-e_j)+D_j\phi _0^-(t+e_k-e_j))
\end{array}
\right )\\
&&=\left (\begin{array}{l}
\zeta _k^{-1}\sum _{j=1}^N
\zeta _jA_j\sum _{t\in\nspace{\widetilde{Z}}{N}_0}\zeta ^tx_0(t)+
\zeta _k^{-1}\sum _{j=1}^N
\zeta _jB_j\sum _{t\in\nspace{\widetilde{Z}}{N}_0}\zeta ^t\phi _0^-(t)\\
\zeta _k^{-1}\sum _{j=1}^N
\zeta _jC_j\sum _{t\in\nspace{\widetilde{Z}}{N}_0}\zeta ^tx_0(t)+
\zeta _k^{-1}\sum _{j=1}^N
\zeta _jD_j\sum _{t\in\nspace{\widetilde{Z}}{N}_0}\zeta ^t\phi _0^-(t)
\end{array}
\right )\\
&&=\zeta _k^{-1}\left (\begin{array}{ll}
\zeta\mathbf{A} & \zeta\mathbf{B}\\
\zeta\mathbf{C} & \zeta\mathbf{D}
\end{array}
\right )\left (\begin{array}{l}
\widehat{x_0}(\zeta )\\
\widehat{\phi _0^-}(\zeta )
\end{array}
\right )=\zeta _k^{-1}\zeta \mathbf{G}\left (\begin{array}{l}
\widehat{x_0}(\zeta )\\
\widehat{\phi _0^-}(\zeta )
\end{array}
\right ).
\end{eqnarray*}
Since the operators $F_\Hspace{X}$ and $F_\Hspace{N^+}$ are unitary, the
operator of multiplication by the block matrix-valued function ``$\cdot
\zeta _k^{-
1}\zeta\mathbf{G}$'':$\mathcal{L}^2(\nspace{T}{N},\Hspace{X})\oplus
\mathcal{L}^2(\nspace{T}{N},\Hspace{N^-
})\to\mathcal{L}^2(\nspace{T}{N},\Hspace{X})\oplus\mathcal{L}^2(\nspace{
T}{N},\Hspace{N^+})$ is contractive. It is easy to show that in this
case for each $\zeta\in\nspace{T}{N}$ both $\zeta _k^{-
1}\zeta\mathbf{G}$ and $\zeta\mathbf{G}$ are contractive operators from
$[\Hspace{X}\oplus\Hspace{N^-},\Hspace{X}\oplus\Hspace{N^+}]$, that
means the dissipativity of $\alpha $, and (iii)$\Rightarrow$(i) holds.
Reversing the argument we obtain (i)$\Rightarrow$(iii).

(i)$\Rightarrow$(ii). As it was shown earlier, (i)$\Rightarrow$(iii),
i.e., from the dissipativity of $\alpha $ the dissipativity of $\alpha _k$
follows for an arbitrary $k$. This means that the generators $W_{\alpha
_k}=W_{\alpha ,k}$ of the associated one-parametric semigroups
are contractive operators for all $k\in\{ 1,\ldots ,N\}$ (see Section
\ref{sec:rem-sys} and Introduction). Hence $\mathfrak{W}_\alpha $ is a
semigroup of
contractions.

(ii)$\Rightarrow$(iii). If $\mathfrak{W}_\alpha $ is a semigroup of
contractions in $\Hspace{H_\alpha }$ then each one-parametric semigroup
$\mathfrak{W}_{\alpha _k}=\{ W_{\alpha ,k}^n\mid n\in\mathbb{Z_+}\}$
is also contractive and hence $\alpha _k$ is a dissipative scattering LSDS
(see Section~\ref{sec:rem-sys}).

The proof of Proposition \ref{prop:dissipative} is complete.
\end{proof}
\begin{rem}\label{rem:phys-dissipative}
Inequalities in \eqref{eq:energ-dis} are multiparametric analogues of
inequalities in
\eqref{eq:energ'} characterizing dissipative scattering LSDS in the case
$N=1$; one may
attach to \eqref{eq:energ-dis} also the physical sense of dissipation of
energy in a
system, and in this case 
 relations for ``powers of incident
and
reflected waves'' and the ``energy'' of inner states of a system are
considered on ``wave fronts'' $|t|=\emph{const}$.
\end{rem}
\begin{thm}\label{thm:tf-n-sys}
The transfer function $\theta _\alpha (z)$ of a dissipative scattering LSDS
$\alpha =(N; \mathbf{A}, \mathbf{B}, \mathbf{C}, \mathbf{D};\Hspace{X},
\Hspace{N^-}, \Hspace{N^+})$
is a holomorphic contractive operator- valued function on
$\nspace{D}{N}$.
\end{thm}
\begin{proof}
Since $L_\mathbf{G}(z)=z\mathbf{G}$ is a contractive
$[\Hspace{X}\oplus\Hspace{N^-},\Hspace{X}\oplus\Hspace{N^+}]$-valued
function on $\nspace{D}{N}$,
$L_\mathbf{A}(z):=z\mathbf{A}=P_\Hspace{X}(z\mathbf{G})|\Hspace{X}$ is a
contractive $[\Hspace{X},\Hspace{X}]$-valued function on
$\nspace{D}{N}$, moreover by virtue of the maximum principle for
holomorphic operator-valued functions (see \cite{S} ) $\| z\mathbf{A}\|
<1$ for $z\in\nspace{D}{N}$. This implies that $\theta _\alpha
(z)=z\mathbf{D}+z\mathbf{C}{(I_\Hspace{X}-z\mathbf{A})}^{-1}z\mathbf{B}$
is well defined and holomorphic on $\nspace{D}{N}$. Fix an arbitrary
$z^0\in\nspace{D}{N}$. Then the one-parametric system $\alpha
_{z^0}=(z^0\mathbf{A},z^0\mathbf{B},z^0\mathbf{C},z^0\mathbf{D};
\Hspace{X}, \Hspace{N^-}, \Hspace{N^+})$ is a dissipative scattering LSDS,
and by Theorem \ref{thm:tf-sys} $\theta _{\alpha _{z^0}}(\lambda
)=\lambda
z^0\mathbf{D}+\lambda z^0\mathbf{C}{(I_\Hspace{X}-\lambda
z^0\mathbf{A})}^{-1}\lambda z^0\mathbf{B}$ is a contractive operator-
valued function on $\mathbb{D}$. The function $\theta _\alpha (z)$ is
holomorphic at $z^0$ and hence continuous at this point. Therefore
\begin{displaymath}
\| \theta _\alpha (z^0)\| =\lim_{\lambda\in\mathbb{D},\ \lambda\to 1}\|
\theta _\alpha (\lambda z^0)\| =\lim_{\lambda\in\mathbb{D},\ \lambda\to
1}\| \theta _{\alpha _{z^0}}(\lambda )\|\leq 1,
\end{displaymath}
that completes the proof.
\end{proof}

\subsection{Multiparametric conservative scattering
LSDS}\label{sec:conserv}
We shall call $\alpha =(N; \mathbf{A}, \mathbf{B}, \mathbf{C},
\mathbf{D};\Hspace{X}, \Hspace{N^-}, \Hspace{N^+})$ a
\emph{multiparametric conservative scattering LSDS} if for any
$\zeta\in\nspace{T}{N}\quad \zeta\mathbf{G}\in
[\Hspace{X}\oplus\Hspace{N^-},\Hspace{X}\oplus\Hspace{N^+}]$
is a unitary operator. It is evident that conservative scattering LSDS
is a special case of dissipative one. For $N=1$ this notion coincides with
the notion of conservative scattering LSDS with
discrete time (see Introduction and Section~\ref{sec:rem-sys}).
\begin{prop}\label{prop:conserv}
Let $\alpha =(N; \mathbf{A}, \mathbf{B}, \mathbf{C},
\mathbf{D};\Hspace{X}, \Hspace{N^-}, \Hspace{N^+})$ be given,
$\mathfrak{W}_\alpha $ and $\alpha _k$ (for some $k\in\{ 1,\ldots ,N\}$)
be the associated semigroup and the associated one-parametric LSDS
respectively. Then the following statements are equivalent:

\emph{(i)} $\alpha$ is a multiparametric conservative scattering LSDS;

\emph{(ii)} $\mathfrak{W}_\alpha $ is a semigroup of unitary operators
in $\Hspace{H}_\alpha $;

\emph{(iii)} $\alpha _k$ is a conservative scattering LSDS;

\emph{(iv)} the $N$-tuple of matrices of the system $\alpha $
\begin{displaymath}
G_k=\left (\begin{array}{ll}
A_k & B_k \\
C_k & D_k
\end{array}
\right )\in [\Hspace{X}\oplus\Hspace{N^-},
\Hspace{X}\oplus\Hspace{N^+}]\quad (k=1, \ldots ,N)
\end{displaymath}
satisfies the following conditions:
\begin{eqnarray}
\sum_{k=1}^NG_k^*G_k &=& I_{\Hspace{X}\oplus\Hspace{N^-
}},\label{eq:cons1} \\
            G_k^*G_j &=& 0, \quad (k\neq j)\label{eq:cons2} \\
\sum_{k=1}^NG_kG_k^* &=&
I_{\Hspace{X}\oplus\Hspace{N^+}},\label{eq:cons3} \\
            G_kG_j^* &=& 0; \quad (k\neq j)\label{eq:cons4}
\end{eqnarray}

\emph{(v)} the spaces $\Hspace{H^+}:=\Hspace{X}\oplus\Hspace{N^+}$ and
$\Hspace{H^-}:=\Hspace{X}\oplus\Hspace{N^-}$ allow decompositions
\begin{equation}\label{eq:decomp}
\Hspace{H^+}=\bigoplus_{k=1}^N\Hspace{H}_k^+, \quad \Hspace{H^-
}=\bigoplus_{k=1}^N\Hspace{H}_k^-
\end{equation}
such that with respect to them matrices $G_k$ have a block structure
\begin{equation}\label{eq:blocks}
(G_k)_{ij}=\left\{ \begin{array}{ll}
G_k^0 & \mbox{if $(i,j)=(k,k)$,}\\
0     & \mbox{if $(i,j)\neq (k,k)$,}
\end{array}
\right. \quad (k=1, \ldots ,N)
\end{equation}
the operator $G^0:=\sum_{k=1}^NG_k$ is unitary and represented by the
block-diagonal matrix $G^0=\emph{diag}(G_1^0,\ldots ,G_N^0)$;

\emph{(vi)} in the assumptions of \emph{(iv)} from
Proposition~\ref{sec:dissipative} for
$\alpha
$, the following equalities are valid for any $n\in\mathbb{Z_+}\setminus
\{ 0\}$:
\begin{equation}\label{eq:energ-cons}
\sum _{|t|=n-1}{\|\phi ^-(t)\| }^2-\sum _{|t|=n}{\|\phi ^+(t)\| }^2=
\sum _{|t|=n}{\| x(t)\| }^2-\sum _{|t|=n-1}{\| x(t)\| }^2
\end{equation}
and in the same assumptions for $\alpha ^*$, for any
$n\in\mathbb{Z_+}\setminus \{ 0\}$:
\begin{equation}\label{eq:energ-cons*}
\sum _{|t|=n-1}{\|\phi _* ^-(t)\| }^2-\sum _{|t|=n}{\|\phi _*^+(t)\|
}^2=
\sum _{|t|=n}{\| x_*(t)\| }^2-\sum _{|t|=n-1}{\| x_*(t)\| }^2
\end{equation}
\end{prop}
\begin{proof}
(i)$\Leftrightarrow$(iv). Rewrite the conservativity conditions for
$\alpha $ as follows:
\begin{eqnarray}
\forall \zeta\in\nspace{T}{N}\quad
{(\zeta\mathbf{G})}^*(\zeta\mathbf{G}) &=&
\sum_{k=1}^N\sum_{j=1}^N\bar{\zeta _k}\zeta
_jG_k^*G_j=I_{\Hspace{X}\oplus\Hspace{N^-}},\label{eq:cons-cond1} \\
\forall \zeta\in\nspace{T}{N}\quad
(\zeta\mathbf{G}){(\zeta\mathbf{G})}^* &=&
\sum_{k=1}^N\sum_{j=1}^N\zeta _k\bar{\zeta
_j}G_kG_j^*=I_{\Hspace{X}\oplus\Hspace{N^+}}.\label{eq:cons-cond2}
\end{eqnarray}
Equating the coefficients under corresponding multipowers of
$\zeta\in\nspace{T}{N}$ (here $\bar{\zeta _k}\zeta _j=\zeta ^{e_j-e_k}$)
we obtain that \eqref{eq:cons-cond1} is equivalent to
\eqref{eq:cons1}$\&$\eqref{eq:cons2}, and \eqref{eq:cons-cond2} is
equivalent to \eqref{eq:cons3}$\&$\eqref{eq:cons4}.

(iv)$\Leftrightarrow$(v). Let (iv) be fulfilled. Denote by
$\overline{\Hspace{Y}}$ the closure of $\Hspace{Y}$. Define
$\Hspace{H}_k^-:=\overline{G_k^*\Hspace{H^+}},\
\Hspace{H}_k^+:=\overline{G_k\Hspace{H^-}}\quad (k=1,\ldots ,N)$. Then
$\Hspace{H}_k^-\perp \Hspace{H}_j^-,\ \Hspace{H}_k^+\perp
\Hspace{H}_j^+$ for $k\neq j$. Indeed, if $h_1^+,\ h_2^+\in\Hspace{H^+}$
then by \eqref{eq:cons4}
$\scalar{G_k^*h_1^+}{G_j^*h_2^+}=\scalar{h_1^+}{G_kG_j^*h_2^+}=0$ for
$k\neq j$, hence by the continuity argument we get $\Hspace{H}_k^-\perp
\Hspace{H}_j^-$ for $k\neq j$. Analogously, $\Hspace{H}_k^+\perp
\Hspace{H}_j^+$ for $k\neq j$. As it was shown earlier,
(iv)$\Leftrightarrow$(i), hence $G^0=\sum_{j=1}^NG_j$ is a unitary
operator. This implies $\Hspace{H^+}=\sum_{j=1}^NG_j\Hspace{H^-}\subset
\sum_{j=1}^N\overline{G_j\Hspace{H^-
}}=\sum_{j=1}^N\Hspace{H}_j^+\subset\Hspace{H^+}$. Thus
$\Hspace{H^+}=\sum_{j=1}^N\Hspace{H}_j^+$ moreover
$\Hspace{H^+}=\bigoplus_{j=1}^N\Hspace{H}_j^+$. Analogously,
$\Hspace{H^-}=\bigoplus_{j=1}^N\Hspace{H}_j^-$, and \eqref{eq:decomp} is
valid.
Since by \eqref{eq:cons4} for any $h^+\in \Hspace{H^+}$ we have
$G_kG_j^*h^+=0$ for
$k\neq j$, by the continuity argument we get $G_k\Hspace{H}_j^-=\{ 0\}$
for $k\neq j$. Then $\Hspace{H}_k^+=\overline{G_k\Hspace{H^-
}}=\overline{G_k\bigoplus_{j=1}^N\Hspace{H}_j^-
}=\overline{G_k\Hspace{H}_k^-}=\overline{\sum_{j=1}^NG_j\Hspace{H}_k^-
}=\overline{G^0\Hspace{H}_k^-}=G^0\Hspace{H}_k^-
=\sum_{j=1}^NG_j\Hspace{H}_k^-=G_k\Hspace{H}_k^-$. Analogously,
$\Hspace{H}_k^-=G_k^*\Hspace{H}_k^+$. Thus, with respect to
\eqref{eq:decomp}
operators $G_k$ have the structure of block matrices \eqref{eq:blocks}
where
$G_k^0=P_{\Hspace{H}_k^+}G_k|\Hspace{H}_k^-\quad (k=1,\ldots ,N)$, and
(iv)$\Rightarrow$(v) is obtained. Conversely, if (v) is fulfilled then
one can verify \eqref{eq:cons1}--\eqref{eq:cons4} immediately, that gives
(v)$\Rightarrow$(iv).

(iii)$\Leftrightarrow$(vi). The proof is analogous to the proof of
(iii)$\Leftrightarrow$(iv) in Proposition \ref{sec:dissipative}, but here it
is necessary
to carry out the same argument twice: both for $\alpha $ and for $\alpha
^*$, and both for the corresponding equalities \eqref{eq:energ-cons} and
\eqref{eq:energ-cons*}.

(iii)$\Leftrightarrow$(i), (i)$\Rightarrow$(ii), (ii)$\Rightarrow$(iii)
can be established similarly to the corresponding parts of
Proposition~\ref{sec:dissipative}.
\end{proof}
\begin{rem}\label{rem:phys-cons}
\eqref{eq:energ-cons} and \eqref{eq:energ-cons*} are multidimensional
analogues of energy equalities
characterizing conservative scattering LSDS in the case $N=1$ and also
mean the conservation of energy under the direct and the inverse
directions of waves
propagation; as for the general case of dissipative system (see Remark
\ref{rem:phys-dissipative})
energy relations are considered here on the ``wave fronts''
$|t|=\emph{const}$.
\end{rem}
\begin{rem}\label{rem:LP-group}
It follows from Proposition \ref{prop:conserv}
 that for a conservative scattering LSDS
$\alpha $ the semigroup $\mathfrak{W}_\alpha =\{ \mathbf{W}_\alpha
^t\mid t\in\nspace{Z}{N}_+\}$ can be extended to the group
$\mathfrak{\widetilde{W}}_\alpha =\{ \mathbf{W}_\alpha ^t\mid
t\in\nspace{Z}{N}\}$ of unitary operators in $\Hspace{H}_\alpha $, and
one can say about the \emph{associated Lax-Phillips group}
$\mathfrak{\widetilde{W}}_\alpha$ for $\alpha $. Indeed, it is
sufficient to show that $W_{\alpha ,k}^{-1}W_{\alpha ,j}=W_{\alpha
,j}W_{\alpha ,k}^{-1}$ for $k\neq j$. But these relations are equivalent
to commutativity relations for $\mathfrak{W}_\alpha $:
$W_{\alpha ,j}W_{\alpha ,k}=W_{\alpha ,k}W_{\alpha ,j}\quad (k\neq j)$
since they follow after the multiplication of the last ones from the
right and from the left by the unitary operator $W_{\alpha ,k}^{-1}$.
\end{rem}

\section{The class of transfer functions of multiparametric
conservative scattering LSDSs}\label{sec:class}
Recall that the \emph{generalized Schur class} $S_N(\Hspace{N^-
},\Hspace{N^+})$ (see \cite{Ag} ) is the class of all holomorphic on
$\nspace{D}{N}$ functions $\theta (z)$ with values in $[\Hspace{N^-
},\Hspace{N^+}]$ where $\Hspace{N^-}$ and $\Hspace{N^+}$ are separable
Hilbert spaces, such that for any separable Hilbert space $\Hspace{Y}$,
for any $N$-tuple $\mathbf{T}=(T_1,\ldots ,T_N)$ of commuting
contractions in $\Hspace{Y}$ and for any positive $r<1$
\begin{eqnarray}
\| \theta (r\mathbf{T})\| & \leq & 1,\label{eq:n-schur}\\
\theta (r\mathbf{T})=\theta (rT_1,\ldots
,rT_N) &:=& \sum_{t\in\nspace{Z}{N}_+}\widehat{\theta _t}\otimes
{(r\mathbf{T})}^t\in [\Hspace{N^-}\otimes \Hspace{Y},\Hspace{N^+}\otimes
\Hspace{Y}],\label{eq:func-cal}
\end{eqnarray}
$\widehat{\theta _t}$ are Maclaurin's coefficients of $\theta (z)$, and
the convergence of series in \eqref{eq:func-cal} is understood in the
sense of norm in
$[\Hspace{N^-}\otimes \Hspace{Y},\Hspace{N^+}\otimes \Hspace{Y}]$. Since
one can choose in particular $\Hspace{Y}=H^2(\nspace{D}{N})$ (the Hardy
space on $\nspace{D}{N}$) and $T_k=$``$\cdot z_k$'' (i.e., the operator
of
multiplication by the $k$-th independent variable in
$\Hspace{Y}=H^2(\nspace{D}{N})$) for $k=1,\ldots ,N$, by virtue of the
isomorphism of the Hilbert spaces $\Hspace{N^\pm}\otimes
H^2(\nspace{D}{N})$ and $H^2(\nspace{D}{N},\Hspace{N^\pm})$ (the last
ones are the Hardy spaces of $\Hspace{N^\pm}$-valued functions on
$\nspace{D}{N}$) according to \eqref{eq:n-schur} we have
\begin{displaymath}
{\| \theta \| }_\infty =\sup_{0<r<1}\|\mbox{``}\cdot \theta
(rz)\mbox{''}\|
_{[H^2(\nspace{D}{N},\Hspace{N^-
}),H^2(\nspace{D}{N},\Hspace{N^+})]}=\sup_{0<r<1}\|\theta
(r\mathbf{T})\| \leq 1
\end{displaymath}
(here ${\| \cdot \| }_\infty$ is the norm in the space $H^\infty
(\nspace{D}{N},[\Hspace{N^-},\Hspace{N^+}])$ of bounded holomorphic on
$\nspace{D}{N}\quad [\Hspace{N^-},\Hspace{N^+}]$-valued functions), thus
$\theta$ turns out to be an element of the unit ball $B_N(\Hspace{N^-
},\Hspace{N^+})$ of the Banach space  $H^\infty
(\nspace{D}{N},[\Hspace{N^-},\Hspace{N^+}])$, and we obtain
\begin{equation}\label{eq:classes}
S_N(\Hspace{N^-},\Hspace{N^+})\subseteq B_N(\Hspace{N^-},\Hspace{N^+}).
\end{equation}
If $N=1$ then due to the von Neumann inequality (see \cite{vN})
$S_N(\Hspace{N^-},\Hspace{N^+})=S(\Hspace{N^-},\Hspace{N^+})$ i.e. the
Schur class. It is known that for $N=1$ (see \cite{vN}) and for $N=2$
(see \cite{An}) one has in fact the sign ``='' in \eqref{eq:classes} for
any
$\Hspace{N^-}$ and $\Hspace{N^+}$, i.e., the classes $S_N(\Hspace{N^-
},\Hspace{N^+})$ and $B_N(\Hspace{N^-},\Hspace{N^+})$ coincide. For
$N>2$, as it follows from \cite{V}, these classes do not coincide for
any $\Hspace{N^-}$ and $\Hspace{N^+}$ different from $\{ 0\}$. As
J.~Agler showed in \cite{Ag} $\theta\in S_N(\Hspace{N^-},\Hspace{N^+})$
if and only if there are separable Hilbert spaces $\Hspace{M}_k$ and
holomorphic functions $F_k:\nspace{D}{N}\to [\Hspace{N^-
},\Hspace{M}_k]\quad (k=1,\ldots ,N)$ such that
\begin{equation}\label{eq:agler}
I_\Hspace{N^-}-{\theta (\lambda )}^*\theta (z)=\sum_{k=1}^N(1-
\overline{\lambda _k}z_k){F_k(\lambda )}^*F_k(z)\quad \mbox{for all
$\lambda ,z\in\nspace{D}{N}$}.
\end{equation}
We shall also use the following result of \cite{Ag}.
\begin{lem}\label{lem:agler}
Let $\Hspace{N},\ \Hspace{K},\ \Hspace{L}$ be separable
Hilbert spaces, $g:\nspace{D}{N}\to [\Hspace{N},\Hspace{K}]$ and
$f:\nspace{D}{N}\to [\Hspace{N},\Hspace{L}]$ be holomorphic functions,
and
\begin{equation}\label{eq:spaces}
\Hspace{G}:=\bigvee_{\lambda\in\nspace{D}{N}}g(\lambda )\Hspace{N},\quad
\Hspace{F}:=\bigvee_{\lambda\in\nspace{D}{N}}f(\lambda )\Hspace{N}.
\end{equation}
If
\begin{equation}\label{eq:g-f}
{g(\lambda )}^*g(z)={f(\lambda)}^*f(z)\quad \mbox{for all $\lambda
,z\in\nspace{D}{N}$}
\end{equation}
then there exists unique isomorphism
$L:\Hspace{G}\to\Hspace{F}$ such that
\begin{equation}\label{eq:L-isom}
f(\lambda )=Lg(\lambda )\quad \mbox{for all $\lambda\in\nspace{D}{N}$.}
\end{equation}
\end{lem}
Denote by $S_N^0(\Hspace{N^-},\Hspace{N^+})$ and
$B_N^0(\Hspace{N^-},\Hspace{N^+})$ the subclasses of $S_N(\Hspace{N^-
},\Hspace{N^+})$ and $B_N(\Hspace{N^-},\Hspace{N^+})$
respectively, consisting of all operator-valued functions vanishing at
$z=0$. By \eqref{eq:tf-n-sys} and
Theorem~\ref{thm:tf-n-sys}, the transfer functions of
multiparametric dissipative scattering LSDSs belong to
$B_N^0(\Hspace{N^-},\Hspace{N^+})$, but we have not a full description
of the transfer functions class for dissipative systems. For conservative
systems we have obtained such description.
\begin{thm}\label{thm:class}
The holomorphic function $\theta :\nspace{D}{N}\to [\Hspace{N^-
},\Hspace{N^+}]$ is the transfer function of some multiparametric
conservative scattering LSDS $\alpha $ (i.e., $\theta =\theta _\alpha $)
if and only if $\theta\in S_N^0(\Hspace{N^-},\Hspace{N^+})$.
\end{thm}
\begin{proof}
Let $\alpha =(N; \mathbf{A}, \mathbf{B}, \mathbf{C},
\mathbf{D};\Hspace{X}, \Hspace{N^-}, \Hspace{N^+})$
be a multiparametric conservative scattering LSDS. Take some $\phi _1^-
,\ \phi _2^-\in \Hspace{N^-}$ and define the input multisequences for
$\alpha $
\begin{displaymath}
\phi _k^-(t)=\left\{ \begin{array}{ll}
\phi _k^- & \mbox{for $t=0$,}\\
0         & \mbox{for $t\in\nspace{\widetilde{Z}}{N}_+\setminus \{
0\}$.}
\end{array}
\right.\quad (k=1,2)
\end{displaymath}
Then by \eqref{eq:n-series}, \eqref{eq:transf-n-state1} and
\eqref{eq:transf-n-output1} the formulas
\begin{displaymath}
\widehat{\phi _k^-}(z)\equiv\phi _k^-,\ \
\widehat{x_k}(z)={(I_\Hspace{X}-
z\mathbf{A})}^{-1}z\mathbf{B}\phi _k^-,\ \ \widehat{\phi _k^+}(z)=\theta
_\alpha (z)\phi _k^-\ \ (k=1,2)
\end{displaymath}
define holomorphic functions on $\nspace{D}{N}$ (note that by
Theorem~\ref{thm:tf-n-sys} $\theta _\alpha $ is a holomorphic on
$\nspace{D}{N}$
operator-valued function). For arbitrary $\lambda ,\ z\in \nspace{D}{N}$
by \eqref{eq:transform-n-sys}--\eqref{eq:transf-n-output1} we get
\begin{eqnarray*}
\scalar{(I_\Hspace{N^-}-{\theta _\alpha (\lambda )}^*\theta
_\alpha (z))\phi _1^-}{\phi _2^-}_\Hspace{N^-}=\scalar{\phi _1^-}{\phi
_2^-}_\Hspace{N^-}-\scalar{\widehat{\phi _1^+}(z)}{\widehat{\phi
_2^+}(\lambda )}_\Hspace{N^+}\\
=\bigscalar{\left (\begin{array}{l}
\widehat{x_1}(z)\\
\phi _1^-
\end{array}
\right )}{\left (\begin{array}{l}
\widehat{x_2}(\lambda )\\
\phi _2^-
\end{array}
\right )}_{\Hspace{X}\oplus\Hspace{N^-}}-\bigscalar{\left
(\begin{array}{l}
\widehat{x_1}(z)\\
\widehat{\phi _1^+}(z)
\end{array}
\right )}{\left (\begin{array}{l}
\widehat{x_2}(\lambda )\\
\widehat{\phi _2^+}(\lambda )
\end{array}
\right )}_{\Hspace{X}\oplus\Hspace{N^+}}\\
=\bigscalar{\left (\begin{array}{l}
\widehat{x_1}(z)\\
\phi _1^-
\end{array}
\right )}{\left (\begin{array}{l}
\widehat{x_2}(\lambda )\\
\phi _2^-
\end{array}
\right )}_{\Hspace{X}\oplus\Hspace{N^-}}-
\bigscalar{z\mathbf{G}\left (\begin{array}{l}
\widehat{x_1}(z)\\
\phi _1^-
\end{array}
\right )}{\lambda\mathbf{G}\left (\begin{array}{l}
\widehat{x_2}(\lambda )\\
\phi _2^-
\end{array}
\right
)}_{\Hspace{X}\oplus\Hspace{N^+}}\\
=\bigscalar{(I_{\Hspace{X}\oplus\Hspace{N^
-}}-{(\lambda\mathbf{G})}^*z\mathbf{G})\left (\begin{array}{l}
\widehat{x_1}(z)\\
\phi _1^-
\end{array}
\right )}{\left (\begin{array}{l}
\widehat{x_2}(\lambda )\\
\phi _2^-
\end{array}
\right )}_{\Hspace{X}\oplus\Hspace{N^-}}.
\end{eqnarray*}
According to Proposition \ref{prop:conserv} (see \eqref{eq:cons1} and
\eqref{eq:cons2}) the last one is
equal to
\begin{eqnarray*}
\bigscalar{(\sum_{k=1}^NG_k^*G_k-\sum_{k=1}^N\overline{\lambda
_k}z_kG_k^*G_k)\left (\begin{array}{l}
\widehat{x_1}(z)\\
\phi _1^-
\end{array}
\right )}{\left (\begin{array}{l}
\widehat{x_2}(\lambda )\\
\phi _2^-
\end{array}
\right )}_{\Hspace{X}\oplus\Hspace{N^-}}=\\
\bigscalar{\sum_{k=1}^N(1-\overline{\lambda _k}z_k)G_k^*G_k
\left (\begin{array}{l}
\widehat{x_1}(z)\\
\phi _1^-
\end{array}
\right )
}{\left (\begin{array}{l}
\widehat{x_2}(\lambda )\\
\phi _2^-
\end{array}
\right
)}_{\Hspace{X}\oplus\Hspace{N^-}}=\\
\bigscalar{\sum_{k=1}^N(1-\overline{\lambda _k}z_k)
{\left (\begin{array}{c}
{(I_\Hspace{X}-\lambda\mathbf{A})}^{-1}\lambda\mathbf{B}\\
I_\Hspace{N^-}
\end{array}
\right )}^*G_k^*G_k
\left (\begin{array}{c}
{(I_\Hspace{X}-z\mathbf{A})}^{-1}z\mathbf{B}\\
I_\Hspace{N^-}
\end{array}
\right )\phi _1^-}{\phi _2^-}_\Hspace{N^-}.
\end{eqnarray*}
Since $\phi _1^-$ and $\phi _2^-$ are arbitrary elements of $\Hspace{N^-
}$ we obtain \eqref{eq:agler} for $\theta =\theta _\alpha $ where
$\Hspace{M}_k:=\Hspace{X}\oplus\Hspace{N^+},\quad F_k(z):=G_k\left
(\begin{array}{c}
{(I_\Hspace{X}-z\mathbf{A})}^{-1}z\mathbf{B}\\
I_\Hspace{N^-}
\end{array}
\right )$, thus $\theta _\alpha \in S_N(\Hspace{N^-},\Hspace{N^+})$.
Moreover, according to \eqref{eq:tf-n-sys} $\theta (0)=0$ and hence
$\theta _\alpha
\in S_N^0(\Hspace{N^-},\Hspace{N^+})$.

Conversely, suppose  $\theta \in S_N^0(\Hspace{N^-},\Hspace{N^+})$.
Then for $\theta $ \eqref{eq:agler} holds with some separable Hilbert
spaces
$\Hspace{M}_k$ and some holomorphic functions $F_k:\nspace{D}{N}\to
[\Hspace{N^-},\Hspace{M}_k]\quad (k=1,\ldots ,N)$. Rewrite
\eqref{eq:agler} in the
following form: $\forall\lambda ,\ z\in\nspace{D}{N}$
\begin{equation}\label{eq:my-g-f}
\sum_{k=1}^N{(\lambda
_kF_k(\lambda ))}^*z_kF_k(z)+I_\Hspace{N^-}=\sum_{k=1}^N{F_k(\lambda
)}^*F_k(z)+{\theta (\lambda )}^*\theta (z).
\end{equation}
Set $\Hspace{N}:=\Hspace{N^-},\quad
\Hspace{M}:=\bigoplus_{k=1}^N\Hspace{M}_k,\quad
\Hspace{K}:=\Hspace{M}\oplus\Hspace{N^-},\quad
\Hspace{L}:=\Hspace{M}\oplus\Hspace{N^+}$, define the functions
$g:\nspace{D}{N}\to [\Hspace{N},\Hspace{K}]$ and $f:\nspace{D}{N}\to
[\Hspace{N},\Hspace{L}]$ by
\begin{displaymath}
g(\lambda ):=\left (\begin{array}{c}
\lambda _1F_1(\lambda )\\
\vdots \\
\lambda _NF_N(\lambda )\\
I_\Hspace{N^-}
\end{array}
\right ),\quad
f(\lambda ):=\left (\begin{array}{c}
F_1(\lambda )\\
\vdots \\
F_N(\lambda )\\
\theta (\lambda )
\end{array}
\right ).\quad
(\lambda \in\nspace{D}{N})
\end{displaymath}
Then for these functions \eqref{eq:my-g-f} means \eqref{eq:g-f}. If we
define the spaces
$\Hspace{G}$ and $\Hspace{F}$ by \eqref{eq:spaces} then we arrive into
the
conditions
of Lemma \ref{lem:agler}, and therefore there is unique unitary operator
$L:\Hspace{G}\to\Hspace{F}$ satisfying \eqref{eq:L-isom}, i.e.,
\begin{equation}\label{eq:my-L-isom}
L\left (\begin{array}{c}
\lambda _1F_1(\lambda )\\
\vdots \\
\lambda _NF_N(\lambda )\\
I_\Hspace{N^-}
\end{array}
\right )=\left (\begin{array}{c}
F_1(\lambda )\\
\vdots \\
F_N(\lambda )\\
\theta (\lambda )
\end{array}
\right )\quad\mbox{for all $\lambda\in\nspace{D}{N}$.}
\end{equation}
Denote by $P_k:=P_{\Hspace{M}_k}$ the orthogonal projector in
$\Hspace{M}$ onto $\Hspace{M}_k$, and set $\mathbf{P}:=(P_1,\ldots
,P_N)$,
\begin{displaymath}
F(\lambda ):=\left (\begin{array}{c}
F_1(\lambda )\\
\vdots \\
F_N(\lambda )
\end{array}
\right )\in [\Hspace{N^-},\Hspace{M}].\quad (\lambda\in\nspace{D}{N})
\end{displaymath}
Then \eqref{eq:my-L-isom} turn into
\begin{equation}\label{eq:my-L-isom1}
L\left (\begin{array}{c}
\lambda\mathbf{P}F(\lambda )\\
I_\Hspace{N^-}
\end{array}
\right )=\left (\begin{array}{c}
F(\lambda )\\
\theta (\lambda )
\end{array}
\right )\quad \mbox{for all $\lambda\in\nspace{D}{N}$.}
\end{equation}
If we substitute $\lambda =0$ into this equality we get
\begin{equation}\label{eq:lambda=0}
\forall\phi ^-\in\Hspace{N^-}\quad L\left (\begin{array}{l}
0\\
\phi ^-
\end{array}
\right )=\left (\begin{array}{c}
F(0)\phi ^-\\
0
\end{array}
\right ).
\end{equation}
>From \eqref{eq:my-L-isom1} and \eqref{eq:lambda=0} we conclude that for
any $\phi ^-\in\Hspace{N^-}$
and for any $\lambda\in\nspace{D}{N}\quad \left (\begin{array}{c}
\lambda\mathbf{P}F(\lambda )\phi ^-\\
0
\end{array}
\right )=\left (\begin{array}{c}
\lambda\mathbf{P}F(\lambda )\\
I_\Hspace{N^-}
\end{array}
\right )\phi ^--\left (\begin{array}{l}
0\\
I_\Hspace{N^-}
\end{array}
\right )\phi ^-\in\Hspace{G}$, and
\begin{equation}\label{eq:lambda-else}
L\left (\begin{array}{c}
\lambda\mathbf{P}F(\lambda )\phi ^-\\
0
\end{array}
\right )=\left (\begin{array}{c}
F(\lambda )-F(0)\\
\theta (\lambda )
\end{array}
\right )\phi ^-.
\end{equation}
Since for any $\phi _1^-,\ \phi _2^-\in\Hspace{N^-}$ and
$\lambda\in\nspace{D}{N}$ vectors $\left (\begin{array}{l}
0\\
\phi _1^-
\end{array}
\right )$ and $\left (\begin{array}{c}
\lambda\mathbf{P}F(\lambda )\phi _2^-\\
0
\end{array}
\right )$ are orthogonal in $\Hspace{G}$, and $L$ is a unitary operator,
by \eqref{eq:lambda=0} and \eqref{eq:lambda-else}  $\left
(\begin{array}{c}
F(0)\phi _1^-\\
0
\end{array}
\right )\perp\left (\begin{array}{c}
F(\lambda )-F(0)\\
\theta (\lambda )
\end{array}
\right )\phi _2^-$ in $\Hspace{F}$, and hence in
$\Hspace{M}\oplus\Hspace{N^+}$. Therefore $F(0)\phi _1^-\perp (F(\lambda
)-F(0))\phi _2^-$ in $\Hspace{M}$. Thus
\begin{equation}\label{eq:span}
\bigvee_{\lambda\in\nspace{D}{N}}F(\lambda )\Hspace{N^-}=\left
(\bigvee_{\lambda\in\nspace{D}{N}}(F(\lambda )-F(0))\Hspace{N^-}\right
)\oplus F(0)\Hspace{N^-}.
\end{equation}
The second summand is a closed lineal in $\Hspace{M}$ because $F(0)$ is
an isometry (this follows from \eqref{eq:agler} for $\lambda =z=0$). Set
$U:=L|\bigvee_{\lambda\in\nspace{D}{N}}\lambda\mathbf{P}F(\lambda
)\Hspace{N^-}$. Then according to \eqref{eq:lambda-else} $U$
is a unitary operator from
$\bigvee_{\lambda\in\nspace{D}{N}}\lambda\mathbf{P}F(\lambda
)\Hspace{N^-}$ onto $\bigvee_{\lambda\in\nspace{D}{N}}\left
(\begin{array}{c}
F(\lambda )-F(0)\\
\theta (\lambda )
\end{array}
\right )\Hspace{N^-}$, and
\begin{equation}\label{eq:u}
U(\lambda\mathbf{P})F(\lambda )=
\left (\begin{array}{c}
F(\lambda )-F(0)\\
\theta (\lambda )
\end{array}
\right )
\quad\mbox{for all $\lambda\in\nspace{D}{N}$.}
\end{equation}
We have the inclusions
$\bigvee_{\lambda\in\nspace{D}{N}}\lambda\mathbf{P}F(\lambda
)\Hspace{N^-}\subset\Hspace{M}$ and
\begin{displaymath}
\bigvee_{\lambda\in\nspace{D}{N}}\left
(\begin{array}{c}
F(\lambda )-F(0)\\
\theta (\lambda )
\end{array}
\right )\Hspace{N^-}\subset (\Hspace{M}\ominus F(0)\Hspace{N^-
})\oplus\Hspace{N^+}.
\end{displaymath}
Enlarging if necessary $\Hspace{M}$ (for
instance, by addition of an infinite countable number to its dimension)
we can achieve
\begin{eqnarray*}
\lefteqn{\mbox{dim}\left (\Hspace{M}\ominus
\bigvee_{\lambda\in\nspace{D}{N}}\lambda\mathbf{P}F(\lambda
)\Hspace{N^-}\right )=}\\
&&\mbox{dim}\left (((\Hspace{M}\ominus
F(0)\Hspace{N^-
})\oplus\Hspace{N^+})\ominus\bigvee_{\lambda\in\nspace{D}{N}}\left
(\begin{array}{c}
F(\lambda )-F(0)\\
\theta (\lambda )
\end{array}
\right )\Hspace{N^-}\right )
\end{eqnarray*}
and extend $U$ to a unitary operator $\widetilde{U}:\Hspace{M}\to
(\Hspace{M}\ominus F(0)\Hspace{N^-})\oplus\Hspace{N^+}$.
Set $\Hspace{X}:=\Hspace{M}\ominus F(0)\Hspace{N^-}$. Then
$\widetilde{U}$ maps $\Hspace{X}\oplus F(0)\Hspace{N^-}$ isometrically
onto $\Hspace{X}\oplus\Hspace{N^+}$, and $\widetilde{U}\left
(\begin{array}{ll}
I_\Hspace{X} & 0 \\
0            & F(0)
\end{array}
\right )$ maps $\Hspace{X}\oplus\Hspace{N^-}$ isometrically onto
$\Hspace{X}\oplus\Hspace{N^+}$. Set
\begin{displaymath}
G_k:=\widetilde{U}P_k\left (\begin{array}{ll}
I_\Hspace{X} & 0 \\
0            & F(0)
\end{array}
\right ).\quad (k=1,\ldots ,N)
\end{displaymath}
As operators from $\Hspace{X}\oplus\Hspace{N^-}$ into
$\Hspace{X}\oplus\Hspace{N^+}$ they have a block form
\begin{displaymath}
G_k=\left
(\begin{array}{ll}
A_k & B_k \\
C_k & D_k
\end{array}
\right ) \quad (k=1,\ldots ,N)
\end{displaymath}
and define the conservative scattering LSDS $\alpha =(N;
\mathbf{A}, \mathbf{B}, \mathbf{C}, \mathbf{D};\Hspace{X}, \Hspace{N^-},
\Hspace{N^+})$. Indeed, for any $\zeta\in\nspace{T}{N}\quad
\zeta\mathbf{G}=\widetilde{U}(\zeta\mathbf{P})\left (\begin{array}{ll}
I_\Hspace{X} & 0 \\
0            & F(0)
\end{array}
\right )$ is a unitary operator from $\Hspace{X}\oplus\Hspace{N^-}$ onto
$\Hspace{X}\oplus\Hspace{N^+}$ because $\left (\begin{array}{ll}
I_\Hspace{X} & 0 \\
0            & F(0)
\end{array}
\right )\in [\Hspace{X}\oplus\Hspace{N^-},\Hspace{M}],\quad
\zeta\mathbf{P}\in [\Hspace{M},\Hspace{M}]$ and $\widetilde{U}\in
[\Hspace{M},\Hspace{X}\oplus\Hspace{N^+}]$ are unitary operators. By
virtue of \eqref{eq:span} and \eqref{eq:u} for an arbitrary
$\lambda\in\nspace{D}{N}$
\begin{eqnarray*}
\lefteqn{\left (\begin{array}{c}
F(\lambda )-F(0)\\
\theta (\lambda )
\end{array}
\right )=U(\lambda\mathbf{P})F(\lambda
)=\widetilde{U}(\lambda\mathbf{P})F(\lambda
)} \\
&&=\widetilde{U}(\lambda\mathbf{P})\left (\begin{array}{c}
F(\lambda )-F(0)\\
F(0)
\end{array}
\right )=\widetilde{U}(\lambda\mathbf{P})
\left (\begin{array}{ll}
I_\Hspace{X} & 0 \\
0            & F(0)
\end{array}
\right )
\left (\begin{array}{c}
F(\lambda )-F(0) \\
I_\Hspace{N^-}
\end{array}
\right ) \\
&&=\lambda\mathbf{G}
\left (\begin{array}{c}
F(\lambda )-F(0) \\
I_\Hspace{N^-}
\end{array}
\right )=
\left (\begin{array}{l}
\lambda\mathbf{A}(F(\lambda ))-F(0))+\lambda\mathbf{B} \\
\lambda\mathbf{C}(F(\lambda ))-F(0))+\lambda\mathbf{D}
\end{array}
\right ).
\end{eqnarray*}
This implies $F(\lambda )-F(0)={(I_\Hspace{X}-\lambda\mathbf{A})}^{-
1}\lambda\mathbf{B}$ and $\theta (\lambda
)=\lambda\mathbf{C}{(I_\Hspace{X}-\lambda\mathbf{A})}^{-
1}\lambda\mathbf{B}+\lambda\mathbf{D}$, i.e., $\theta =\theta _\alpha $,
that completes the proof.
\end{proof}

Consider for $\alpha =(N; \mathbf{A}, \mathbf{B}, \mathbf{C},
\mathbf{D};\Hspace{X}, \Hspace{N^-}, \Hspace{N^+})$ the following
subspace
\begin{equation}\label{eq:cc-n-sys}
\Hspace{X}_1:=\bigvee_{p,\ k,\
j}p(\mathbf{A},\mathbf{A}^*)(B_k\Hspace{N^-
}+C_j^*\Hspace{N^+})
\end{equation}
of $\Hspace{X}$,
where $p$ runs the set of all monomials in $2N$ non-commuting variables,
$k$ and $j$ run the set $\{ 1,\ldots ,N\}$ (cf. \eqref{eq:cc-sys}). We
shall call
LSDS $\alpha $ \emph{closely connected} if $\Hspace{X}_1=\Hspace{X}$.
\begin{thm}\label{thm:cc-n}
Let $\alpha =(N; \mathbf{A}, \mathbf{B}, \mathbf{C},
\mathbf{D};\Hspace{X}, \Hspace{N^-}, \Hspace{N^+})$ be an arbitrary
multiparametric LSDS. Then:

\emph{(i)} the decomposition $\Hspace{X}=\Hspace{X}_0\oplus\Hspace{X}_1$
takes place where $\Hspace{X}_1$ is defined by \eqref{eq:cc-n-sys}; with
respect to
this decomposition
\begin{equation}\label{eq:cc-matr}
\begin{array}{ll}
A_k=\left (\begin{array}{ll}
A_k^{(0)} & 0 \\
0         & A_k^{(1)}
\end{array}
\right ), & B_k=\left (\begin{array}{l}
0 \\
B_k^{(1)}
\end{array}
\right ), \\
C_k=(\begin{array}{ll}
0 & C_k^{(1)}
\end{array}), & D_k=D_k^{(1)},
\end{array}
\quad (k=1,\ldots ,N)
\end{equation}
and $\alpha _1:=(N; \mathbf{A}^{(1)},
\mathbf{B}^{(1)}, \mathbf{C}^{(1)}, \mathbf{D}^{(1)};\Hspace{X}_1,
\Hspace{N^-}, \Hspace{N^+})$ is a closely connected multiparametric
LSDS;

\emph{(ii)} if $\alpha $ is a dissipative (resp. conservative) scattering
LSDS then $\alpha _1$ is also a dissipative (resp. conservative) scattering
LSDS;

\emph{(iii)} $\theta _{\alpha _1}=\theta _\alpha $.
\end{thm}
\begin{proof}
It is clear from \eqref{eq:cc-n-sys} that $\Hspace{X}_1$ is the minimal
subspace in
$\Hspace{X}$ containing $B_k\Hspace{N^-}$ and $C_k^*\Hspace{N^+}$, and
reducing $A_k$ for all $k\in\{ 1,\ldots ,N\}$. This implies
\eqref{eq:cc-matr}, and
(i) follows. If $\alpha $ is a dissipative system then for any
$\zeta\in\nspace{T}{N}$
\begin{displaymath}
\zeta\mathbf{G}^{(1)}=\left (\begin{array}{ll}
\zeta\mathbf{A}^{(1)} & \zeta\mathbf{B}^{(1)} \\
\zeta\mathbf{C}^{(1)} & \zeta\mathbf{D}^{(1)}
\end{array}
\right
)=P_{\Hspace{X}_1\oplus\Hspace{N^+}}(\zeta\mathbf{G})|\Hspace{X}_1\oplus
\Hspace{N^-}
\end{displaymath}
is a contractive operator, thus $\alpha _1$ is also dissipative.
If $\alpha $ is a conservative system then for any
$\zeta\in\nspace{T}{N}\quad \zeta\mathbf{G}\in
[\Hspace{X}\oplus\Hspace{N^-},\Hspace{X}\oplus\Hspace{N^+}]$ is a
unitary operator, and we get from \eqref{eq:cc-matr}:
\begin{equation}\label{eq:cc-pencil}
\zeta\mathbf{G}=\left (\begin{array}{ll}
\zeta\mathbf{A}^{(0)} & 0 \\
0                     & \zeta\mathbf{G}^{(1)}
\end{array}
\right )\in [\Hspace{X}_0\oplus (\Hspace{X}_1\oplus\Hspace{N^-
}),\Hspace{X}_0\oplus (\Hspace{X}_1\oplus\Hspace{N^+})].
\end{equation}
This implies that $\zeta\mathbf{G}^{(1)}$ is a unitary operator for any
$\zeta\in\nspace{T}{N}$, thus $\alpha _1$ is a conservative system, and
the statement (ii) follows. Further, for $z\in\nspace{D}{N}\quad
(z\mathbf{C}^{(1)})(z\mathbf{B}^{(1)})=(z\mathbf{C})P_{\Hspace{X}_1}(z
\mathbf{B})=(z\mathbf{C})(z\mathbf{B})$ because $z\mathbf{B}\Hspace{N^-
}\subset\Hspace{X}_1$. For $z\in\nspace{D}{N}$ and for
$n\in\mathbb{Z_+}\quad
z\mathbf{C}^{(1)}{(z\mathbf{A}^{(1)})}^nz\mathbf{B}^{(1)}=z\mathbf{C}P_{
\Hspace{X}_1}{(z\mathbf{A}P_{\Hspace{X}_1})}^nz\mathbf{B}=z\mathbf{C}{(z
\mathbf{A})}^nz\mathbf{B}$ because $z\mathbf{B}\Hspace{N^-
}\subset\Hspace{X}_1$ and $z\mathbf{A}\Hspace{X}_1\subset\Hspace{X}_1$.
Therefore $\theta _{\alpha
_1}(z)=z\mathbf{D}^{(1)}+z\mathbf{C}^{(1)}{(I_\Hspace{X}-
z\mathbf{A}^{(1)})}^{-
1}z\mathbf{B}^{(1)}=z\mathbf{D}+z\mathbf{C}{(I_\Hspace{X}-
z\mathbf{A})}^{-1}z\mathbf{B}=\theta _\alpha (z)$. The proof of
Theorem~\ref{thm:cc-n} is complete.
\end{proof}
\begin{rem}\label{rem:cc-main}
If $\alpha $ is conservative then due to \eqref{eq:cc-pencil} the
restriction of the linear pencil $L_\mathbf{A}(\zeta
)=\zeta\mathbf{A}\in
[\Hspace{X},\Hspace{X}]\quad (\zeta\in\nspace{T}{N})$ on
$\Hspace{X}_0\subset\Hspace{X}$, that is the pencil
$L_{\mathbf{A}^{(0)}}(\zeta )=\zeta\mathbf{A}^{(0)}\in
[\Hspace{X}_0,\Hspace{X}_0]\quad (\zeta\in\nspace{T}{N})$,
consists of unitary operators.
\end{rem}
We shall call the linear pencil $L_\mathbf{A}(\zeta )\quad
(\zeta\in\nspace{T}{N})$ of contractive operators in a Hilbert space
$\Hspace{X}$ \emph{completely non-unitary} if there is no any proper
subspace $\Hspace{X}_0$ in $\Hspace{X}$ reducing $L_\mathbf{A}(\zeta )$
for each $\zeta\in\nspace{T}{N}$ (or, equivalently, reducing all $A_k$
for $k=1,\ldots ,N$) and such that a pencil $L_{\mathbf{A}^{(0)}}(\zeta
)=L_\mathbf{A}(\zeta )|\Hspace{X}_0\quad (\zeta\in\nspace{T}{N})$ is
unitary (i.e. consisting of unitary operators).
\begin{thm}\label{thm:c-non-unitary}
The conservative scattering LSDS $\alpha =(N; \mathbf{A}, \mathbf{B},
\mathbf{C}, \mathbf{D};
\Hspace{X},\\ \Hspace{N^-},
 \Hspace{N^+})$ is
closely connected if and only if the pencil $L_\mathbf{A}(\zeta )\quad
(\zeta\in\nspace{T}{N})$ is completely non-unitary.
\end{thm}
\begin{proof}
Suppose that $\alpha $ is closely connected and there is a proper
subspace $\Hspace{X}_0\subset\Hspace{X}$ reducing $L_\mathbf{A}(\zeta )$
for all $\zeta\in\nspace{T}{N}$ so that $L_{\mathbf{A}^{(0)}}(\zeta
)=L_\mathbf{A}(\zeta )|\Hspace{X}_0\quad (\zeta\in\nspace{T}{N})$ is a
unitary pencil of operators. Then
$\Hspace{X}=\Hspace{X}_0\oplus\Hspace{X}_0^\perp$
and, with respect to this decomposition, \eqref{eq:cc-matr} and
\eqref{eq:cc-pencil} take place.
But in this case $\Hspace{X}_0^\perp =\Hspace{X}_1\neq\Hspace{X}$ where
$\Hspace{X}_1$ is defined by \eqref{eq:cc-n-sys}. So we get the
contradiction with
the supposition on the close connectedness of $\alpha $. Thus
$L_\mathbf{A}(\zeta )\quad (\zeta\in\nspace{T}{N})$ is a completely 
non-unitary pencil of operators.

Conversely, let $L_\mathbf{A}(\zeta )\quad (\zeta\in\nspace{T}{N})$ be a
completely non-unitary pencil of operators. Then by
Remark~\ref{rem:cc-main} in the
decomposition
$\Hspace{X}=\Hspace{X}_0\oplus\Hspace{X}_1\quad\Hspace{X}_0=\{ 0\}$ is
necessary, i.e., $\Hspace{X}_1=\Hspace{X}$ and $\alpha $ is closely
connected.
\end{proof}
Theorem~\ref{thm:class} together with Theorem~\ref{thm:cc-n} imply the
following
multidimensional analogue of Theorem~\ref{thm:cc-cons}.
\begin{thm}\label{thm:n-cc-cons}
An arbitrary operator-valued function $\theta\in S_N^0(\Hspace{N^-
},\Hspace{N^+})$ allows a closely connected conservative realization,
i.e., there is such a multiparametric closely connected conservative
scattering LSDS $\alpha $ that $\theta =\theta _\alpha $.
\end{thm}
\begin{rem}
Theorem~\ref{thm:n-cc-cons} is only partial generalization of
Theorem~\ref{thm:cc-cons} to the case
$N>1$, since Theorem~\ref{thm:cc-cons} claims also uniqueness of a
closely connected
conservative realization for $\theta\in S^0(\Hspace{N^-},\Hspace{N^+})$
up to unitary similarity, that means in particular that the state
spaces of two such realizations are isomorphic.
\end{rem}
The following example shows that for $N>1$ there are closely connected
conservative realizations of one and the same $\theta\in
S_N^0(\Hspace{N^-},\Hspace{N^+})$ with non-isomorphic state spaces.
\begin{example}\label{ex}
Define $\theta (z)=z_1z_2$ for $z\in\nspace{D}{2}$. Evidently, $\theta
(z)$ is a scalar contractive holomorphic function on $\nspace{D}{2}$ and
$\theta (0)=0$, i.e., $\theta\in
B_2^0(\mathbb{C},\mathbb{C})=S_2^0(\mathbb{C},\mathbb{C})$. Define
$\alpha
:=(2;\mathbf{A},\mathbf{B},\mathbf{C},\mathbf{D};\mathbb{C},\mathbb{C},
\mathbb{C})$ where $A_1=A_2=0,\quad B_1=0,\quad B_2=1,\quad C_1=1,\quad
C_2=0,\quad D_1=D_2=0$, i.e., for $z=(z_1,\ z_2)\in\nspace{C}{2}$
\begin{displaymath}
z\mathbf{G}=\left (\begin{array}{ll}
z\mathbf{A} & z\mathbf{B} \\
z\mathbf{C} & z\mathbf{D}
\end{array}
\right )=\left (\begin{array}{ll}
0   & z_2 \\
z_1 & 0
\end{array}
\right ).
\end{displaymath}
It is clear that if $z\in\nspace{T}{2}$ then $z\mathbf{G}$ is a unitary
operator in $\nspace{C}{2}$, i.e., $\alpha $ is a conservative system.
Since $B_2\Hspace{N^-}=\mathbb{C}$, from relations
$\mathbb{C}=\Hspace{X}\supset\Hspace{X}_1\supset B_2\Hspace{N^-
}=\mathbb{C}$ we obtain the close connectedness of $\alpha $. Finally,
it is clear that $\theta _\alpha (z)=\theta (z)=z_1z_2$.

Now define $\alpha
'=(2;\mathbf{A'},\mathbf{B'},\mathbf{C'},\mathbf{D'};\nspace{C}{3},
\mathbb{C},\mathbb{C})$ where
\begin{displaymath}
\begin{array}{ll}
A_1'=\left (\begin{array}{ccc}
0 & 0          & -1/\sqrt{2} \\
0 & 0          & 0 \\
0 & 1/\sqrt{2} & 0
\end{array}
\right ), & B_1'=\left (\begin{array}{c}
1/\sqrt{2} \\
0 \\
0
\end{array}
\right ),\\
C_1'=(\begin{array}{ccc}
0 & 1/\sqrt{2} & 0
\end{array}
), & D_1'=0, \\
A_2'=\left (\begin{array}{ccc}
0 & 0          & 0 \\
0 & 0          & 1/\sqrt{2} \\
-1/\sqrt{2}&0  & 0
\end{array}
\right ), & B_2'=\left (\begin{array}{c}
0 \\
1/\sqrt{2} \\
0
\end{array}
\right ), \\
C_2'=(\begin{array}{ccc}
1/\sqrt{2} & 0 & 0
\end{array}
), & D_2'=0,
\end{array}
\end{displaymath}
i.e., for $z=(z_1,\ z_2)\in\nspace{C}{2}$
\begin{displaymath}
z\mathbf{G'}=\left (\begin{array}{cc}
z\mathbf{A'} & z\mathbf{B'} \\
z\mathbf{C'} & z\mathbf{D'}
\end{array}
\right )=\left (\begin{array}{cccc}
0             & 0            & -z_1/\sqrt{2} & z_1/\sqrt{2} \\
0             & 0            & z_2/\sqrt{2}  & z_2/\sqrt{2} \\
-z_2/\sqrt{2} & z_1/\sqrt{2} & 0             & 0  \vspace{1mm}     \\
z_2/\sqrt{2}  & z_1/\sqrt{2} & 0             & 0
\end{array}
\right ).
\begin{picture}(10,10)
\put(-200,-11){\line(1,0){180}}
\put(-58,-27){\line(0,1){60}}
\end{picture}
\end{displaymath}
Evidently, if $z\in\nspace{T}{2}$ then $z\mathbf{G'}$ is a unitary
operator in $\nspace{C}{4}$, i.e., $\alpha '$ is a conservative system.
Since $B_1'\Hspace{N^-}=\mathbb{C}\oplus\{ 0\}\oplus\{
0\}\subset\nspace{C}{3},\quad B_2'\Hspace{N^-}=\{
0\}\oplus\mathbb{C}\oplus\{ 0\}\subset\nspace{C}{3},\quad
A_2'B_1'\Hspace{N^-}=\{ 0\}\oplus\{
0\}\oplus\mathbb{C}\subset\nspace{C}{3}$, from relations
$\nspace{C}{3}=\Hspace{X}\supset\Hspace{X}_1\supset B_1'\Hspace{N^-
}+B_2'\Hspace{N^-}+A_2'B_1'\Hspace{N^-}=\nspace{C}{3}$ we obtain the
close connectedness of $\alpha '$. Since for $z\in \nspace{C}{2}$ we
have $(z\mathbf{C'})(z\mathbf{B'})=z_1z_2$, and for any
$z\in\nspace{C}{2}$ and $n\in\mathbb{Z_+}\setminus \{ 0\}$ we have
${(z\mathbf{A'})}^nz\mathbf{B'}=0$, we obtain $\theta _{\alpha
'}(z)=z\mathbf{D'}+z\mathbf{C'}{(I_{\nspace{C}{3}}-z\mathbf{A'})}^{-
1}z\mathbf{B'}=(z\mathbf{C'})(z\mathbf{B'})=z_1z_2=\theta (z)$. Thus
both $\alpha $ and $\alpha '$ are closely connected conservative
realizations of $\theta (z)$, however
$1=\mbox{dim}\Hspace{X}\neq\mbox{dim}\Hspace{X'}=3$.
\end{example}
\begin{rem}\label{rem:lack-uniq}
The lack of uniqueness of a closely connected conservative realization
for a given function from the generalized Schur class
 is not a specific character
of our approach: a similar phenomenon appears in \cite{BT}.
\end{rem}
\begin{rem}\label{rem:cc-min}
The closely connected conservative realization $\alpha =(N; \mathbf{A},
\mathbf{B},\\
 \mathbf{C}, \mathbf{D};\Hspace{X}, \Hspace{N^-},
\Hspace{N^+})$ of $\theta\in S_N^0(\Hspace{N^-},\Hspace{N^+})$ is a
minimal conservative realization in the following sense: if
$\Hspace{X'}$ is some proper subspace in $\Hspace{X},\quad
A_k'=P_\Hspace{X'}A_k|\Hspace{X'},\quad B_k'=P_\Hspace{X'}B_k,\quad
C_k'=C_k|\Hspace{X'},\quad D_k'=D_k$ for $k=1,\ldots ,N$ then the system
$\alpha '=(N; \mathbf{A'}, \mathbf{B'}, \mathbf{C'},
\mathbf{D'};\Hspace{X'}, \Hspace{N^-}, \Hspace{N^+})$ can not be
conservative, not to mention a conservative realization of $\theta$.
Indeed, in the opposite case the following representation takes place:
\begin{displaymath}
\zeta\mathbf{G}=\left (\begin{array}{ll}
\zeta\mathbf{A}^{(0)} & 0 \\
0                     & \zeta\mathbf{G'}
\end{array}
\right )\quad (\zeta\in\nspace{T}{N})
\end{displaymath}
where $\zeta\mathbf{A}^{(0)}\in
[\Hspace{X}\ominus\Hspace{X'},\Hspace{X}\ominus\Hspace{X'}]$ and
$\zeta\mathbf{G'}\in [\Hspace{X'}\oplus\Hspace{N^-
},\Hspace{X'}\oplus\Hspace{N^+}]\quad (\zeta\in\nspace{T}{N})$ are
linear pencils of unitary operators, and we get
$\Hspace{X}_1\subset\Hspace{X'}\neq\Hspace{X}$ (see \eqref{eq:cc-n-sys})
that contradicts to the close connectedness of $\alpha$.
\end{rem}

\emph{Acknowledgements.} I wish to express my gratitude to D.Z.~Arov for
inspiring discussions and his continual interest to this my work, and to
N.~Young who told me about J.~Agler's paper and sent it to me.

\emph{Research supported in part by INTAS Grant 93-0322-ext and by
Ukrainian-Israeli project of scientific co-operation
(contract no. 2M/1516-97).}

\bibliographystyle{amsplain}
\bibliography{msys}

\providecommand{\bysame}{\leavevmode\hbox to3em{\hrulefill}\thinspace}
\begin{thebibliography}{10}

\bibitem{AA}
V.M. Adamjan and D.Z. Arov, \emph{On unitary couplings of semiunitary
  operators}, Mat. Issled., Kishinev \textbf{1} (1966), no.~2, 3--66,
  (Russian).

\bibitem{Ag}
J.~Agler, \emph{On the representation of certain holomorphic functions defined
  on a polydisc}, Topics in Operator Theory: {Ernst D. Hellinger} Memorial
  Volume (L.~de~Branges, I.~Gohberg, and J.~Rovnyak, eds.), Oper. Theory and
  Appl., vol.~48, Birkh\"{a}user-Verlag, Basel, 1990, pp.~47--66.

\bibitem{An}
T.~Ando, \emph{On a pair of commutative contractions}, Acta Sci. Math. (Szeged)
  \textbf{24} (1963), 88--90.

\bibitem{AFK}
M.A. Arbib, P.L. Falb, and R.E. Kalman, \emph{Topics in mathematical system
  theory}, {McGraw}-Hill, New York, 1969.

\bibitem{A1}
D.Z. Arov, \emph{Scattering theory with dissipation of energy}, Dokl. Akad.
  Nauk \textbf{216} (1974), no.~4, 716--718, (Russian).

\bibitem{A2}
\bysame, \emph{Passive linear stationary dynamic systems}, Sibirsk. Math. Zh.
  \textbf{20} (1979), no.~2, 211--228, (Russian).

\bibitem{BC}
J.A. Ball and N.~Cohen, \emph{De {Branges-Rovnyak} operator models and systems
  theory: a survey}, Topics in Matrix and Operator Theory (H.~Bart, I.~Gohberg,
  and M.A. Kaashoek, eds.), Oper. Theory Adv. Appl., vol.~50,
  Birkh\"{a}user-Verlag, Basel, 1991, pp.~93--136.

\bibitem{BT}
J.A. Ball and T.Trent, \emph{Unitary colligations, reproducing kernel {Hilbert}
  spaces, and {Nevanlinna-Pick} interpolation in several variables}, Preprint.

\bibitem{Be}
M.F. Bessmertniy, \emph{Functions of several variables in the theory of finite
  linear structures}, Candidate dissertation, Kharkov University, Kharkov,
  1982, (Russian).

\bibitem{Bo}
N.K. Bose, \emph{Problems and progress in multidimensional systems theory},
  Proc. {IEEE} \textbf{65} (1977), no.~6, 824--840.

\bibitem{Br}
M.S. Brodskii, \emph{Unitary operator colligations and their characteristic
  functions}, Uspekhi Mat. Nauk \textbf{33} (1978), no.~4, 141--168, (Russian).

\bibitem{CS}
M.~Cotlar and C.~Sadosky, \emph{Integral representations of bounded {Hankel}
  forms defined in scattering systems with a multiparametric evolution group},
  Oper. Theory Adv. Appl. \textbf{35} (1988), 357--375.

\bibitem{GS}
D.~Ga\c{s}par and N.~Suciu, \emph{On the structure of isometric semigroups},
  Oper. Theory Adv. Appl. \textbf{14} (1984), 125--139.

\bibitem{HM}
D.J. Hill and P.J. Moylan, \emph{Dissipative dynamical systems: Basic
  input-output and state properties}, J. Franklin Inst. \textbf{309} (1980),
  327--357.

\bibitem{K}
T.~Kaczorek, \emph{Two-dimensional linear systems}, {LNCIS}, vol.~68,
  Springer-Verlag, Berlin, 1985.

\bibitem{LP}
P.~Lax and R.~Phillips, \emph{Scattering theory}, Academic Press, New York,
  1967.

\bibitem{LI}
M.S. Liv\v{s}ic, \emph{Operator waves in {Hilbert} spaces and related partial
  differential equations}, Integral Equations Operator Theory \textbf{2}
  (1979), no.~1, 25--47.

\bibitem{LKMV}
M.S. Liv\v{s}ic, N.~Kravitsky, A.S. Markus, and V.~Vinnikov, \emph{Theory of
  commuting nonselfadjoint operators}, Mathematics and Its Applications, vol.
  332, Kluwer, Dordrecht, 1995.

\bibitem{LW}
M.S. Liv\v{s}ic and L.~Waksman, \emph{Commuting nonselfadjoint operators in
  {Hilbert} space}, Lect. Notes in Math., vol. 1272, Springer-Verlag, New York,
  1987.

\bibitem{N}
M.A. Nudelman, \emph{The sufficient conditions of the absolute stability of
  optimal passive scattering systems}, Algebra i Analiz \textbf{6} (1994),
  no.~4, 187--203, (Russian).

\bibitem{R}
W.~Rudin, \emph{Function theory in the unit ball of $\nspace{C}{n}$},
  Springer-Verlag, Berlin-Heidelberg-New York, 1980.

\bibitem{S}
L.~Schwartz, \emph{Analyse math\'{e}matique}, vol.~II, Hermann, Paris, 1967.

\bibitem{Se}
M.A. {Semyonov-Tyan-Shanskii}, \emph{Harmonic analysis on {Riemanian} symmetric
  surfaces of a negative curvature and scattering theory}, Izv. Akad. Nauk
  {USSR} Ser. Math. \textbf{40} (1976), no.~3, 562--592, (Russian).

\bibitem{Su}
I.~Suciu, \emph{On the semi-groups of isometries}, Stud. Math. \textbf{30}
  (1968), 101--110.

\bibitem{Sz.-NF}
B.~{Sz.-Nagy} and C.~Foia\c{s}, \emph{Harmonic analysis of operators on
  {Hilbert} spaces}, North Holland, Amsterdam, 1970.

\bibitem{V}
N.~Varopoulos, \emph{On an inequality of von {Neumann} and an application of
  the metric theory of tensor products to operator theory}, J. Funct. Anal.
  \textbf{16} (1974), 83--100.

\bibitem{vN}
J.~von Neumann, \emph{Eine {Spectraltheorie} f\"{u}r allgemeine {Operatoren}
  eines unit\"{a}ren {Raumes}}, Math. Nachr. \textbf{4} (1951), 258--281.

\bibitem{W}
J.C. Willems, \emph{Dissipative dynamical systems}, Archive for Rational
  Mechanics and Analysis \textbf{45} (1972), 321--393.

\bibitem{Z}
V.A. Zolotarev, \emph{Lax-{Phillips} scattering scheme on groups and a
  functional model of {Lie} algebra}, Mat. Sbornik \textbf{183} (1992), no.~5,
  115--144, (Russian).

\end{thebibliography}
$\begin{array}{c}
\emph{DMITRIY S. KALYUZHNIY} \\
\emph{Department of Higher Mathematics} \\
\emph{Odessa State Academy of Civil Engineering} \\
\emph{and Architecture} \\
\emph{Didrihson str. 4} \\
\emph{270029, Odessa} \\
\emph{UKRAINE}
\end{array}$

\end{document}